\documentclass[11pt]{amsart}

\usepackage[margin=3cm]{geometry}

\usepackage[foot]{amsaddr}
\usepackage{amsmath}
\usepackage{amssymb}
\usepackage{amsthm}
\usepackage{graphicx}
\usepackage{float}
\usepackage{color}
\usepackage{ulem}
\usepackage{dsfont}
\usepackage{epstopdf}
\usepackage{hyperref}

\newtheorem{theorem}{Theorem}[section]

\newtheorem{lemma}[theorem]{Lemma}

\newtheorem{conjecture}{Conjecture}

\theoremstyle{definition}

\newcommand{\ep}{\varepsilon}

\title[Selected topics on  reaction-diffusion-advection models]
      {Selected topics on  reaction-diffusion-advection models from spatial 
      ecology}

\thanks{{K.-Y. Lam, Y. Lou}: Department of Mathematics, Ohio State University, Columbus, OH 43210, USA}

\thanks{{S. Liu}: Institute for Mathematical Sciences, Renmin University of China, Beijing 100872, China.}

\author[King-Yeung Lam et al.]
{King-Yeung Lam, \ \ Shuang Liu, \ \ Yuan Lou}


 
\begin{document}
\maketitle

\begin{abstract}
We discuss the effects of movement and spatial heterogeneity on population dynamics via reaction-diffusion-advection models, focusing on the persistence, competition, and evolution of organisms in spatially heterogeneous environments. Topics include Lokta-Volterra competition models, river models, evolution of biased movement, phytoplankton growth, and spatial spread of epidemic disease. Open problems and conjectures are presented. Parts of this survey are adopted from the materials in  \cite{LamLou2019, Lou2015, Lou2019}, and some very recent progress are also included.
\end{abstract}


\section{\bf Introduction}\label{S1}
Recent years have witnessed unprecedented progress of  experimental technologies in the life sciences.
The explosion of empirical  results   have rapidly generated massive sets of loosely structured
data. The analytical methods from mathematics and statistics are required  to synthesize the large data sets and
extract insightful information from them. 
This trend continues to accelerate the development in mathematical biology.  
The study of mathematical biology usually includes two aspects. On the one hand,  by introducing and analyzing  mathematical models, researchers can elucidate and predict the basic mechanisms of underlying biological processes. On the other hand,  deeper understanding of these mechanisms may 
drive the discovery of
new mathematical problems, new analytical techniques, and even launch new  research directions.
Most, if not all, areas of mathematics have found applications in mathematical biology.

As an important mathematical area, nonlinear partial differential equations have  been one of the most active research fields in 
the  21st century. It has also found new opportunities in mathematical biology in recent years.
In \cite{F2010} 
Dr. Avner Friedman discussed some challenging problems arising from mathematical biology, and 
one of the major mathematical tools he used  is nonlinear partial differential equations. 
The {subject} of partial differential equations concerns the changes of
quantities in space and time.
For biology, the importance of space is hardly a question. 
With mathematics, we can endeavor to quantify the effect of space in different biological scenarios.
In this survey, we consider some issues in spatial ecology via reaction-diffusion models,
focusing upon the effect of movement of species on population dynamics  in spatially heterogeneous environments.  An important ongoing trend  
in population dynamics is the  integration with other directions of biosciences, such as  evolution, epidemiology, cell biology, cancer research,  etc.  Our main goal is 
to showcase the role of reaction-diffusion models by integrating across different research directions in mathematical biology, including spatial ecology, evolution and disease transmissions.
We will illustrate that the quantitative analysis of some important questions in biosciences does bring  interesting and novel mathematical problems, which in turn calls for the development of 
new mathematical tools. 
We will also formulate a number of potential research questions along the way.

This paper is organized as follows. In Section \ref{S2} we discuss several  single species models, for which the dynamics have important implications
for the invasion of exotic species. Section \ref{S3} is devoted to the study of  two-species competition models, with emphasis on the effects of dispersal and spatial heterogeneity on the outcome of competition for two similar species. In Section \ref{S4} we investigate the competition models with directed movement and show how different dispersal (either random or biased) can bring dramatic changes  to  population dynamics.  Section \ref{S5} concerns the continuous trait models and evolution of dispersal.
Finally, in Sections \ref{S6} and \ref{S7} we discuss some 
 progress on dynamics of phytoplankton growth and disease transmissions, respectively.

\section{\bf Single species models}\label{S2}

The study of mathematical models for single species has 
{played a primary role}
in  population dynamics. It 
{is also crucial}
in analyzing  the dynamics of multiple  interacting species, 
 e.g. on issues concerning  the invasions of  exotic  species.
In this section we  will focus on two 
types of 
single species models 
 and study 
the effect of varying the diffusion rate. 

\subsection{Logistic model}\label{S2.1}
In this subsection, we consider the 
logistic 
model 
 \begin{equation}\label{eq:1}
\left\{
\aligned
&\partial_t u =d \Delta u + u [m(x)-u] && (x,t)\in \Omega \times (0,
\infty),\\
&\frac{\partial u}{\partial \nu} =0
&&(x,t)\in
\partial\Omega \times (0, \infty),\\
&u(x,0)=u_0(x) &&x\in\Omega,
\endaligned
\right.
\end{equation}
 where $\Omega$ is a bounded domain in $\mathbb{R}^N$ with smooth boundary $\partial\Omega$ and $\nu(x)$ denotes the
unit outward normal vector at $x\in\partial\Omega$. Here $u(x,t)$ represents the density of the single species at location $x$ and time $t$. Parameter $d>0$ is the diffusion rate and $\Delta=\sum_{i=1}^N \frac{\partial^2}{\partial x_i^2}$ is the usual Laplace operator. The function $m\in C^2({\overline\Omega})$ accounts for the local carrying capacity or the intrinsic growth rate of the species, which is assumed to be  strictly positive in ${\overline\Omega}$ for the sake of clarity. The Neumann boundary condition 
indicates that no individuals can move  across
$\partial\Omega$, i.e. the 
 habitat is closed.
We assume that the initial data $u_0$ is non-negative
and not identically zero. 

It is  well known
\cite{CC2003} that 
$$\lim_{t\rightarrow\infty }u(x,t)=u^*(x,d) \quad \text{ uniformly for } x\in{\overline\Omega},$$
where for each $d>0$, the function $u^*(\cdot,d)$ is the unique positive solution of 
\begin{equation}\label{eq:2}
\left\{
\aligned
&d\Delta u+ u(m-u)=0 &&x\in\Omega,
\\
&\frac{\partial u}{\partial \nu}=0 &&x\in\partial\Omega.
\endaligned
\right.
\end{equation}
 If the spatial environment is homogeneous, i.e. $m$ is a positive constant, then 
 $u^*(x,d)\equiv m$ is the unique positive solution of 
 \eqref{eq:2}, which is independent of the diffusion rate $d$. If the spatial environment is heterogeneous,
 i.e. $m=m(x)$
is a non-constant function, then  $u^*(x,d)$ depends non-trivially on  $d$. In this case, a natural question is how the 
total biomass of the species at equilibrium, i.e. $\int_\Omega u^*(x,d)\, \mathrm{d}x$, depends on the diffusion rate $d$. 

The following  property  was 
 observed in
\cite{lou2006}:
\begin{lemma} 
\label{lemma:1}
 If $m$ is non-constant, then for any diffusion rate  $d>0$, 
$$
\int_\Omega u^*>\int_\Omega m=\lim_{d\to 0} \int_\Omega u^*
=\lim_{d\to \infty} \int_\Omega u^*.
$$
\end{lemma}

 Lemma \ref{lemma:1} means that {\rm(i)} the heterogeneous environment
can support a total 
biomass greater 
than the total carrying capacity of
the environment, which is quite different from the homogeneous case;  {\rm(ii)} the  total biomass,
as a function of
$d$, is non-monotone; it  is maximized at some intermediate diffusion rate,  and is
minimized at $d=0$ and $d=\infty$,  respectively. Examples were constructed in \cite{LiangLou2012}  to illustrate that this function can 
{have two or more} 
local maxima. 
 

Lemma \ref{lemma:1} implies that $\int_\Omega u^*/\int_\Omega m>1$ for all $d>0$, and the lower bound $1$ is optimal. For the upper bound of $\int_\Omega u^*/\int_\Omega m$,  Ni raised the following question:
\begin{conjecture}[\cite{Ni2011-a}]\label{conjecture1}
There exists some constant $C=C(N)$ depending only on $N$ such that  
$$
\int_\Omega u^* \Big{/} \int_\Omega m\leq C(N), \quad \text{ and }\quad C(1)=3.$$ 
\end{conjecture}

Conjecture \ref{conjecture1} suggests that  the supremum of the ratio of the  biomass 
to  the  total carrying capacity 
depends only on the spatial dimension. For one-dimensional case, $C(1)=3$   was proved in \cite{Bai2014} and it turns out to be  optimal. However, it is recently shown by  Inoue and Kuto \cite{IK2020} that Conjecture \ref{conjecture1} fails for the higher-dimensional case even when 
$\Omega$ is a multi-dimensional ball.
 That is, the nonlinear operator which maps $m$ to $u^*$
is not bounded as a map from 
the positive cone of $L^1(\Omega)$
to itself.

A related question is the optimal control of total biomass \cite{Ding2010}:


\medskip
\noindent{\bf Optimization problem.} Fix any $\delta\in (0, 1)$, and define the control set 
$$\begin{array}{c}
   U:=\left\{m\in L^{\infty}(\Omega):0\le m\le 1,\,\int_\Omega m=\delta|\Omega|\right\}.
\end{array}
$$ 
 Determine those $m^*\in U$ which can maximize the total 
 biomass $\int_\Omega u^*$ in $U$. 
\medskip


 Nagahara and Yanagida \cite{Nagahara}  showed that, under a regularity assumption, the optimal $m^*$ is of ``bang-bang'' type, i.e. $m^*=\chi_E$  for some measurable set $E\subset\Omega$, where $\chi_E$ denotes the characteristic function of $E$. More recently, Mazari et al.  \cite{MNP2020} proved that the bang-bang property holds for all large enough diffusion rates.
Biologically, the set $E$ can be viewed as the protected area, and the
condition $\int_\Omega m=\delta|\Omega|$ means that the total resources are limited.
The results in \cite{MNP2020, Nagahara} imply that under the limited  resources, the way to maximize the total population size of species is to place all the resources evenly in
some suitable subset $E$ of $\Omega$. However,  the characterization of the set $E$ is a challenging problem. It seems that the problem has not been completely solved even in the one-dimensional case; see
\cite{KLY2008} for some  progress on cylindrical domains.
Recent work \cite{NLY}
for the discrete patch model  suggests that the
set $E$ is periodically fragmented; see also \cite{MR2020} for general fragmentation phenomenon.
For the biased movement model with dispersal term $\nabla\cdot((1+m)\nabla u)$, the optimal distribution of resources for maximizing the survival ability of a species was considered in \cite{MNP2020-2}. 
They showed that the problem has  ``regular" solutions only when the domain is a
ball and the optimal distribution can be characterized in this case; see \cite{CDP2017} for the one-dimensional case.

Compared to \eqref{eq:2}, a more general model is 
\begin{equation}\label{eq:1a}
\left\{
\aligned
&d\Delta u+ r(x) u\left(1-\frac{u}{K(x)}\right)=0 \quad &&x\in\Omega,
\\
&\frac{\partial u}{\partial\nu}=0\quad &&x\in\partial\Omega,
\endaligned
\right.
\end{equation}
where $r$ is the intrinsic growth rate and $K$ is the carrying capacity of environment.
Is it possible to find general sufficient conditions for $r,K$ such that the maximal
biomass size of species can be reached at some 
intermediate  diffusion rate?  This question
for \eqref{eq:1a} was  proposed and studied 
by  DeAngelis  et al. \cite{DeAngelis2015}.
We refer to \cite{HLLN} for some relevant developments in this regard, which  shows that the  total
biomass 
can be a monotone increasing function of the diffusion rate for some choices of $r$ and $K$.

In the study of a predator-prey model \cite{Lou2017}, 
the following problem arises: Is the maximum  of the density in \eqref{eq:2}, i.e.  $\max_{x\in{\overline\Omega}} u^*(x,d)$,  monotone decreasing   with respect to  the diffusion rate $d$? Some recent progress is obtained in \cite{LiRui}, where the monotonicity was proved for several classes of function $m$. However, it remains open to show the monotonicity for 
general function $m$. 
In contrast, the minimum  of the density, i.e.  $\min_{x\in{\overline\Omega}} u^*(x,d)$, is not necessarily monotone increasing in $d$; see \cite{HeXiaoqing-NiWeiMing2016}. We 
conjecture that the difference between the maximum and minimum values of the density, which measures 
the spatial variations of the density, is a monotone decreasing function of $d$.
It is known that 
$\int_\Omega |\nabla u^*|^2$, which also measures 
the spatial variations of the density, is monotone decreasing with respect to $d$; see \cite{LiangLou2012} for more details. 

\subsection{Single species models in rivers}\label{S2.2}

{How do} populations persist in streams when they are constantly subject to  downstream drift?
This problem is 
{intensified when the habitat quality} at the downstream end is very poor. Once the species are washed downstream, the chances of survival will be greatly reduced. This  problem, termed as the ``drift paradox”, has received considerable attention \cite{HPPK1993,M1954}. Speirs and Gurney \cite{Speirs} pointed out that the action of 
dispersal can permit persistence in an advective environment, and they proposed the following reaction-diffusion-advection model: 
\begin{equation}\label{eq:3}
\left\{
\aligned
\smallskip
&\partial_t u=d\partial_{xx} u-\alpha\partial_x u+u(r-u) \,\, &&x\in (0,L), t>0,
\\ \smallskip
&d\partial_x u(0,t)-\alpha u(0,t)=u(L,t)=0
&&t>0,
\\
&u(x,0)=u_0(x) &&x\in (0,L),
\endaligned
\right.
\end{equation}
where positive constants $\alpha$ and $r$ are the advection rate and the growth rate of species, respectively.  The interval $[0,L]$ represents an idealized river, assuming that the upstream end $x=0$ is a no-flux or closed boundary, while 
{at the downstream end $x=L$ the zero Dirichlet boundary condition (also called lethal boundary condition) is imposed. Biologically interpreted, the downstream end describes} the junction of rivers and seas, as freshwater fishes  cannot survive in the seawater.
The persistence of species is equivalent to 
{the instability of the trivial equilibrium, which is in turn characterized by the negativity of the principal eigenvalue of the following problem:} 
 \begin{equation}\label{eq:3.1}
 \left\{
 \medskip
\aligned
&d\partial_{xx}\varphi-\alpha\partial_{x}\varphi+(r+\lambda)\varphi=0\quad x\in(0,L),
\\
& d\partial_x\varphi(0)-\alpha\varphi(0)=\varphi(L)=0.
\endaligned
\right.
 \end{equation}

\begin{figure}
  \centering
  \includegraphics[width=0.6\linewidth]{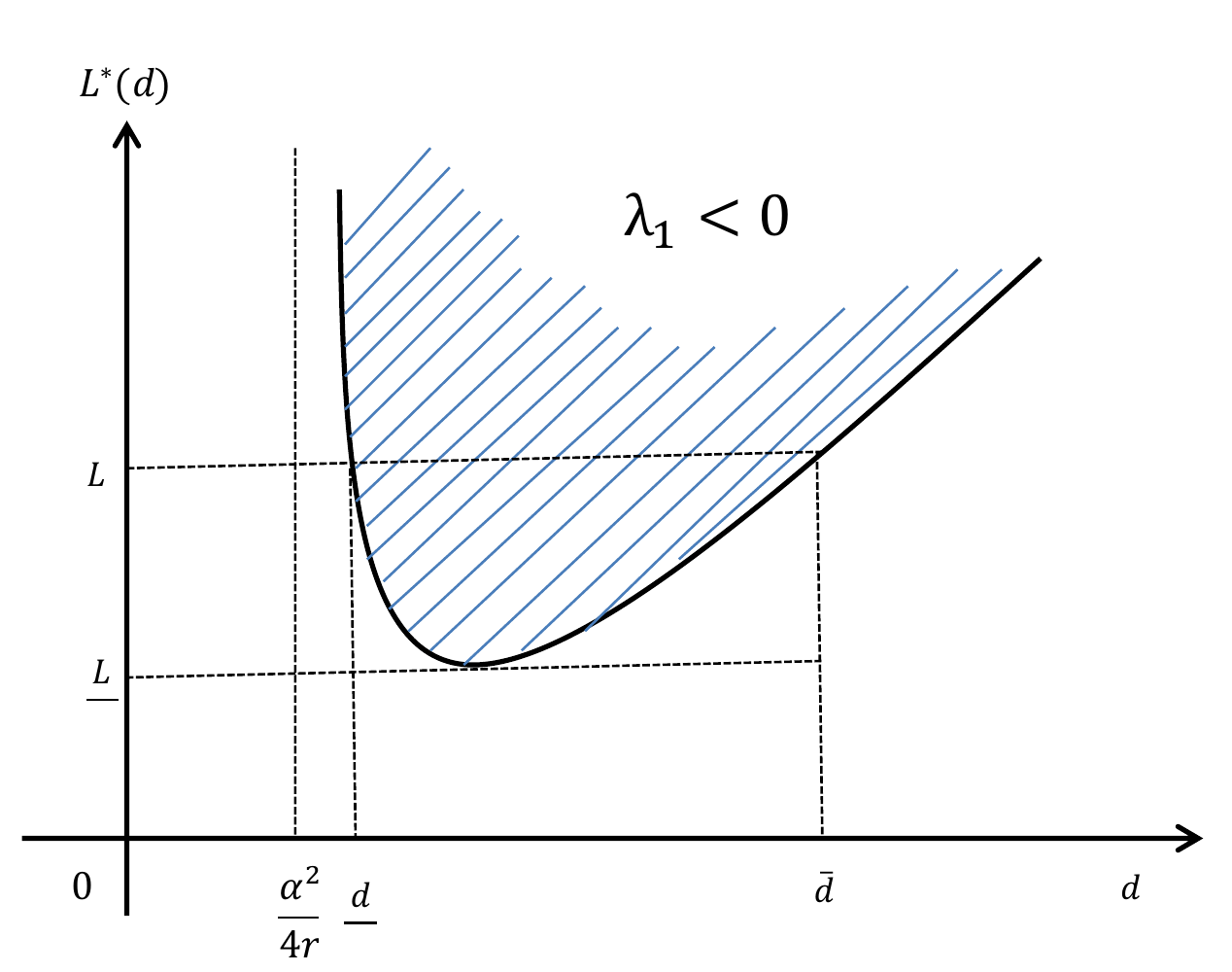}
 \caption{The diagram of  $L^{*}(d)$ illustrated by the black colored solid curve. 
 The parameter region where the species persists is shaded in blue. 
}\label{liufig1}
 \end{figure}

By an explicit calculation,  
Speirs and  Gurney \cite{Speirs} obtained the following necessary and sufficient condition for the persistence of species: 
 \begin{equation}\label{condition}
   \begin{split}
d>\frac{\alpha^2}{4r}
\quad\text{and}\quad
L>L^{*}(d):=\frac{\pi-\arctan\left(\sqrt{4d r-\alpha^2}/\alpha\right)}{\sqrt{4d r-\alpha^2}/2d}.
    \end{split}
   \end{equation}
Set $\underline{L}:=\inf_{d>0}L^{*}(d)$. As illustrated in   Figure \ref{liufig1}, condition \eqref{condition} implies that
 \begin{itemize}
    \item[{\rm(i)}]  if $L\leq \underline{L}$, then the species cannot persist  for any diffusion rate $d$. In this case, the equilibrium $u=0$  is globally stable among all non-negative initial values;
    \smallskip
   \item[{\rm(ii)}] if $L>\underline{L}$, then there exist $0<\underline{d}<\overline{d}$ such that the species persists if and only if $d\in (\underline{d},\overline{d})$. Furthermore, 
   if $d\not\in (\underline{d},\overline{d})$, then $u=0$  is globally stable; if $d\in (\underline{d},\overline{d})$， then \eqref{eq:3} admits a unique positive steady state which is globally stable.
    \end{itemize}

Biologically, the conditions for species survival are that the river has a minimal length and the diffusion rate  of  the species lies within a certain range: when the diffusion rate is small, the drift is so dominant that the species 
is washed downstream  and cannot survive;  when the diffusion rate is large, the species
is exposed too much to the lethal boundary 
at the downstream end and will 
also go to extinction. 
Therefore, it was conjectured in \cite{LamLou2019, LouLutscher2014} that  there is an intermediate diffusion rate  $d^*$ which is the ``best"  strategy, i.e. once the diffusion rate of a species is $d^*$,  no other species can replace it.
This  conjecture implies that if the majority of individuals in a population adopt the strategy $d^*$, 
{then any invasion fails so long as the invading species is adopting a strategy different from $d^*$.}
In the evolutionary game theory, such $d^*$ is referred to as an evolutionary stable strategy \cite{ESS}.
See Subsection \ref{S3.3} and Conjecture \ref{conjecture2} 
for further details.

In  \eqref{eq:3},  the downstream end  is assumed to be  the lethal boundary. Recent work considering general boundary conditions \cite{LouZhou2015, Vasilyeva2011} showed that the species persists if and only if the velocity of flow drift is  less than 
some threshold value.  Under some other boundary conditions, this threshold value turns out to be a monotone decreasing function of the diffusion rate, while this is not the case for \eqref{eq:3}. 
 This leads to a natural evolutionary question: Is the fast diffusion  or slow  diffusion more beneficial to the species in streams?
 We will discuss this question in Subsection \ref{S3.3}. 
 Some new persistence criteria has been found for general boundary conditions at the downstream in a most recent work \cite{HLL2020},
where the critical river length can be either monotone decreasing in the diffusion rate $d$ or it can first decrease and then increase as $d$ increases,  depending on the loss rate at the downstream; see also \cite{XWT2020} for related discussions.  The diagram in Figure \ref{liufig1} can be regarded as  an extreme case of the new  persistence criteria  found in \cite{HLL2020}.  

Finally,  we refer to \cite{DLPZ2020, JPS2019, 
Lutscher2005, Lutscher2006,  Vasilyeva20019, WS2019, WSW2019} for related work on the persistence problem in rivers and  river networks.

\section{\bf Competing species 
models}\label{S3}

One basic question in spatial ecology is: How does the spatial heterogeneity of environment affect the invasion of exotic species and the competition outcome of native and exotic species? 
The term heterogeneity refers to the non-uniform distribution of various environmental conditions. For example, in the oceans, light intensity decreases with respect to the depth
due to light absorption by phytoplankton and water \cite{K1994},  and thus is  unevenly distributed.
An exotic species can successfully invade often because it has  some competitive advantage over the resident species. 
 However, in the  heterogeneous environment, with the help of diffusion alone, invasion is still possible even though the exotic species does not 
{have any apparent}
 competitive advantage.  
 In this section, 
 we consider some classical competition models to illustrate this phenomenon.

\subsection{Lotka-Volterra competition models}\label{S3.1}
We first consider the following   Lotka-Volterra competition model in the homogeneous environment:
\begin{equation}\label{eq:3-liu}
\left\{
\aligned
&\partial_t u=d\Delta u+ u(a_1-b_1u-c_1v) \quad &&(x,t)\in\Omega\times (0, \infty),
\\
&\partial_t u=D\Delta v+ v(a_2-b_2u-c_2v) \quad &&(x,t)\in\Omega\times (0, \infty),
\\
&\frac{\partial u}{\partial\nu}=\frac{\partial v}{\partial\nu}=0\quad &&(x,t)\in\partial\Omega\times (0, \infty),
\\
&u(x,0)=u_0(x), \ v(x, 0)=v_0(x) \quad &&x\in\Omega,
\endaligned
\right.
\end{equation}
where $u(x, t)$ and $v(x,t)$ are the population densities of two competing species 
with  diffusion rates  $d$ and $D>0$, parameters 
$a_1, a_2$ are their intrinsic growth rates, $b_1,c_2$ are the intra-specific competition coefficients, and $b_2, c_1$ are the inter-specific competition coefficients.  
We assume for the moment that $a_i, b_i, c_i$, $i=1,2$, are positive constants, i.e. the spatial environment  is  homogeneous.
If the coefficients  satisfy the weak competition condition 
\begin{equation}\label{eq:3-1}
\frac{b_1}{b_2}>\frac{a_1}{a_2}>\frac{c_1}{c_2},
\end{equation}
then \eqref{eq:3} admits a unique positive steady state $(u^*, v^*)$ given by 
\begin{equation}\label{eq:3-2}
u^*=\frac{a_1c_2-a_2c_1}{b_1c_2-b_2c_1}, \quad v^*=\frac{a_2b_1-a_1b_2}{b_1c_2-b_2c_1}.
\end{equation}
\begin{theorem}[{{\rm \cite{Brown1980}}}]\label{law:1-liu}
Suppose that the environment  is  homogeneous and \eqref{eq:3-1} holds. Then the positive steady state $(u^*, v^*)$ is globally asymptotically stable among all non-negative and non-trivial initial conditions.
\end{theorem}

Theorem \ref{law:1-liu} shows that when the  environment is homogeneous and the competition  is weak,  two competing  species can always coexist regardless of the initial conditions. A natural question arises: Does Theorem \ref{law:1-liu} remain valid in spatially heterogeneous environment?

Obviously, when $a_i, b_i, c_i$ are  non-constant functions in $\Omega$, $(u^*, v^*)$ given by  \eqref{eq:3-2} is generally {non-constant}, and thus not necessarily  a positive steady state of \eqref{eq:3-liu}. 
However we may still ask: If  condition \eqref{eq:3-1} holds pointwisely on $\Omega$,  does  \eqref{eq:3-liu} admit a  unique positive steady state which is globally asymptotically stable? Surprisingly, even if \eqref{eq:3-1} holds pointwisely in  $\Omega$,  it is possible that one of the two competing species  eventually becomes extinct, regardless of  initial conditions, provided that two  competitors have proper diffusion rates.
Mathematically, for the  homogeneous case, \eqref{eq:3-liu} can be viewed as an ordinary differential equation, which has a globally stable positive steady state for each $x\in\Omega$. In contrast,  for the  heterogeneous case, the global attractor 
of the partial differential equation model \eqref{eq:3-liu}  could simply be a semi-trivial steady state, so that the density of one competitor stays bounded from below while the other tends to zero, as $t\to\infty$.
Thus, the spatial heterogeneity of environment also plays a crucial role here. 

For further  discussions on the competitive exclusion in spatially heterogeneous environment, we consider the following simplified competition model:
\begin{equation}\label{eq:4}
\left\{
\aligned
&\partial_t u=d\Delta u+ u{[m(x)-u-bv]} \quad &&(x,t)\in\Omega\times (0, \infty),
\\
&\partial_t v=D\Delta v+ v{[m(x)-cu-v]} \quad &&(x,t)\in\Omega\times (0, \infty),
\\
&\frac{\partial u}{\partial\nu}=\frac{\partial v}{\partial\nu}=0\quad &&(x,t)\in\partial\Omega\times (0, \infty),
\\
&u(x,0)=u_0(x), \ v(x, 0)=v_0(x) \quad &&x\in\Omega,
\endaligned
\right.
\end{equation}
where $m(x)$ represents the common resource of two competitors and  the weak competition condition \eqref{eq:3-1} is reduced  to $0<b, c<1$.

We assume that $m(x)$ is positive, smooth and non-constant to reflect the inhomogeneous distribution of resources across the habitat. Then for $d,D>0$, \eqref{eq:4} has two semi-trivial states, $(u^*, 0)$ and  $(0, v^*)$, where $u^*$ is the unique positive solution of \eqref{eq:2},
while $v^*$ is the unique positive solution of \eqref{eq:2} with $d$ replaced by $D$.  For  the spatially heterogeneous environment, the stability of semi-trivial states 
{is a subtle question}.
For example, the stability of $(u^*, 0)$ is determined by the sign of the principal eigenvalue $\lambda_1$ of the problem
\begin{equation}\label{eq:5}
\left\{
\aligned
&D\Delta \varphi+ (m-cu^*)\varphi+\lambda\varphi=0 \quad &&x\in\Omega,
\\
& \frac{\partial \varphi}{\partial\nu}=0\quad &&x\in\partial\Omega,
\endaligned
\right.
\end{equation}  
{where $u^*(x,d)$ depends implicitly on $d$.}
It is well-known \cite{CC2003} that if  $\lambda_1>0$,  then $(u^*, 0)$ is stable; if $\lambda_1<0$, then $(u^*, 0)$ is unstable. 
When $m$ is constant, under weak competition condition $0 < b,c<1$, direct calculations give $\lambda_1<0$, i.e. $(u^*, 0)$ is always unstable for any $d,D>0$. In contrast, for non-constant $m$, the situation becomes complicated and interesting: First, the variational characterization of the principal eigenvalue  \cite{CC2003} implies   that $\lambda_1$ is monotone increasing with respect to $D$, and it tends to  $\min_{{\overline\Omega}}(cu^*-m)$ as $D$ approaches $0$. Notice from \eqref{eq:2} that $m-u^*$ must change sign in $\Omega$. If $c<1$ (in fact $c\le 1$ is sufficient), then $\min_{{\overline\Omega}}(cu^*-m)<0$, which implies that $\lambda_1<0$ for small $D$.  
For large $D$ we have
$$\lim_{D\to\infty} \lambda_1=\frac{1}{|\Omega|}\left(c\int_\Omega u^*-\int_\Omega m\right),$$
which 
{highlights the role played by} the total 
biomass $\int_\Omega u^*$ 
in the stability of $(u^*,0)$.
To be more precise, define
$$
c_*:=\inf_{d>0} \frac{\int_\Omega m}{\int_\Omega u^*}.
$$
From Subsection \ref{S2.1} we see that $c_*\in (0, 1)$, and if $c<c_*$,  then $(u^*, 0)$ is unstable for all $d,D>0$. It is interesting that when $c\in (c_*, 1)$, we can find some non-empty open set $\Sigma$ such that the semi-trivial solution $(u^*, 0)$ is stable if and only of $(d, D)\in \Sigma$, so that we can define
\begin{equation}\label{sigma}
   \Sigma:=\left\{(d,D)\in (0,\infty)\times (0,\infty): (u^*, 0) \,\text{is stable}\right\}. 
\end{equation}
Observe that for any $ c\in(c_*, 1)$, the set $\Sigma$ is always a subset of $\{(d, D): 0 <d <D\}$.

Recently, He and Ni \cite{HeNi3}  proved a surprising and insightful result: for a class of competition models, if the semi-trivial or positive solution is locally stable, then it must be globally stable. The following result, as a corollary of their general result, completely solves a conjecture proposed  in \cite{lou2006, Lou2007}:

\begin{theorem}\label{slow-dispersal} 
Suppose that $0<b\leq 1$, $c\in (c_*, 1)$ and $m$ is non-constant. If $(d, D)\in \Sigma$, 
then $(u^*, 0)$ is globally stable among all non-negative and non-trivial initial conditions; 
if $(d, D)\not\in \overline{\Sigma}$ and $d<D$, then \eqref{eq:4} admits a unique positive steady state which is globally stable.
\end{theorem}

Theorem \ref{slow-dispersal} shows that for some diffusion rates,  
{competition exclusion takes place} even if the weak competition condition is satisfied everywhere in $\Omega$.  The case $b=1$ and $c_*<c<1$ should  be pointed out particularly. In this case, when the environment is homogeneous, i.e. $m$ is constant,  the equilibrium solution $(0, m)$ is globally stable, so that the {faster diffusing} species 
{excludes the slower one}.  
In contrast, for the spatially heterogeneous environment, Theorem \ref{slow-dispersal} implies that for proper diffusion coefficients,  the final outcome of 
the competition is reversed: 
{the slower diffusing species excludes the faster one!}

Furthermore, when $b=1$ 
and $c$ approaches $1$ from below, the following result holds:
\begin{theorem}\label{slow-dispersal-1}
Assume that $b=c=1$.
If $d<D$, then $(0, v^*)$ is unstable and $(u^*, 0)$ is globally stable. Conversely, if $d>D$, then $(0, v^*)$ is globally stable.
\end{theorem} 

{In fact}, as  $c$ tends to 1 from below, the set  $\Sigma$ increases monotonically and converges to the set 
$\{(d, D): 0<d<D\}$ {in the Hausdorff sense}. Biologically, this implies that the slower diffusing species will replace the faster one in the spatially heterogeneous environment. This is the so-called 
evolution of slow dispersal. 
More precisely, if we consider the slower diffuser  as a mutant among the faster diffusing population, then mutations for slower diffusion rate are selected,
 i.e.  the diffusion rate of the  ``new resident" is smaller than the ``previous resident". Through constant mutation and competition, the diffusion of new species  thus becomes slower and slower. This phenomenon was first discovered by Hastings \cite{Hastings1983}, and Theorem \ref{slow-dispersal-1} was proved by Dockery et al. \cite{DHMP}. Moreover, it is shown in  \cite{Hutson2001} that this  result may fail for spatio-temporally varying environments; see also \cite{LLS2020}.

The above discussion illustrates that {diffusion in a spatially heterogeneous environment has a dramatic}  
effect on the dynamics of competition models, which is rather different from {the case of a homogeneous environment}.
It should be noted that there are still many questions
which remain open even for two competing species. We refer the interested readers to \cite{CS2020, HeNi1, HeNi2, LamNi2012, Ni2011} for some related discussions. For general resource distributions, it is challenging  to fully determine  the stability of the semi-trivial state of \eqref{eq:3-liu}; see \cite{CCL2004, HLMP, 
LMP2006} for examples of multiple reversals of the competition outcomes for {two-species} Lotka-Volterra competition models. 
A remaining problem is to find some conditions on the resource distribution, so that the local stability of two semi-trivial steady states of \eqref{eq:3-liu} can be accurately characterized.
The answer may depend upon how to give some sufficient conditions ensuring that the total 
biomass of a single  species as a function of  diffusion rate has a unique local maximum (and thus also a global maximum).
Another interesting question is what happens if all the coefficients $a_i, b_i, c_i$ in \eqref{eq:3-liu} are non-constant, and  the weak competition condition \eqref{eq:3-1} still holds for all  $x\in \Omega$. See the recent progress in 
\cite{Ni2020} via construction of special Lyapunov functions. 

For other population models, such as the predator-prey models, we can enquire similar issues: How do  diffusion and spatial heterogeneity affect the dynamics of predator-prey models? For the competition models, under the weak competition condition, if the diffusion rate of a species lies within some proper ranges,  the remaining average resources can be negative, which is closely related to the phenomenon that the total biomass of a single species can be greater than the total amount of resources, as discussed in Subsection \ref{S2.1}. 
 Therefore, the fast diffusing invader cannot invade successfully in this case as  their available resources can be  negative. 
 Similar consideration is applicable for the predator-prey models. 
However, the study of predator-prey models turns out to be more complicated \cite{Lou2017}:
Understanding the relationship between the total biomass of species and the total amount of resources is far from being sufficient, and other issues are also  involved, such as the dependence of the maximum density of  species on the diffusion rate \cite{LiRui}.

The study of reaction-diffusion models not only uncover some interesting biological phenomena, but also bring new and interesting mathematical problems. In the next two subsections, we continue on the theme
of competition, but focus on the effect of biased movement of species in the spatially heterogeneous  environment.

\subsection{Competition in rivers}\label{S3.3} 
In this subsection we consider the following two-species competition model in a stream:
\begin{equation}\label{eq:6}
\left\{
\aligned
&\partial_tu=d \partial_{xx}u -\alpha_1 \partial_xu + u (r-u-v) && (x,t)\in (0,L)\times (0, \infty),
\\
 &\partial_tv = D\partial_{xx}v -\alpha_2 \partial_xv + v (r-u-v) && (x,t)\in (0,L)\times (0, \infty),
\\
& d \partial_xu-\alpha_1 u= D \partial_xv-\alpha_2 v=0 \quad &&(x,t)\in\{0,L\}\times (0, \infty), \\
& u(x,0)=u_0(x), \ v(x, 0)=v_0(x)
&& x\in(0, L),
\endaligned
\right.
\end{equation}
where $\alpha_1$ and $\alpha_2$ represent the advection rates  of species $u$ and $v$, respectively. When 
$\alpha_1=\alpha_2=0$, it is shown in Subsection \ref{S3.1} that the species with the smaller  diffusion rate will drive the
other species to extinction, provided that $r$ is non-constant.  What happens if $\alpha_1,\alpha_2\neq0$?

For the case when $r$ is a positive constant (i.e. the environment is homogeneous) and $\alpha_1=\alpha_2=\alpha>0$,
the following result was established in \cite{LouZhou2015}: 
\begin{theorem}
\label{law:2}
Suppose that $r>0$ is constant and $\alpha>0$, $D>d$. Then the semi-trivial state $(0, v^*)$ is globally stable, where $v^*$ is the unique positive solution of
\begin{equation}\label{eq:7}
\left\{
\aligned
&D\partial_{xx}v -\alpha \partial_xv + v (r-v)=0  &&  x\in (0,L),
\\
&D\partial_xv-\alpha v=0  &&   x=0, L.
\endaligned
\right.
\end{equation}
\end{theorem}

Theorem \ref{law:2} shows that under the influence of  passive drift, 
the species with faster diffusion will be favored, which is contrary to Theorem \ref{slow-dispersal-1} for model \eqref{eq:4} in the non-advective environment.   In the  river, the ``fate" of the two competitors is reversed: the slower diffusing species  will  be eventually wiped out! An intuitive reasoning is that the passive drift will push more individuals downstream to force a concentration of population  at the downstream end, which causes mismatch of the population density with the resource distribution. Faster diffusion can counterbalance such mismatch by sending more individuals upstream and therefore, it is 
competitively advantageous for species to adopt faster diffusion.

When $r$ is constant, the case  $\alpha_1\neq\alpha_2$ and $D=d$ was considered in \cite{LouXiaoZhou},
where it was proved that the smaller advection is always 
beneficial. 
A similar explanation is that under the  drift, the species will move towards the downstream end, leading to the overcrowding and the overmatching of resources there, which in turn  results in the extinction of the species with the stronger advection  that can  send more individuals to the downstream. The general case $D\neq d$ and $\alpha_1\neq\alpha_2$ was addressed in \cite{Zhou2016},  which showed that  the fast diffusion and weak advection  will be favored, consistent with  the findings in \cite{LouXiaoZhou, LouZhou2015}. 
 
 When $r$ is non-constant, i.e. the spatial environment is heterogeneous, 
the dynamics of \eqref{eq:6} turns out to be more complicated.
The simplest situation might be 
when $r$ is decreasing, for which we have
the following conjecture:

\begin{conjecture}\label{conjecture2-a} 
Suppose that $r$ is positive and decreasing in $(0, L)$. Then 
there exists some $\alpha^*>0$, independent of $d$ and $D$,  such that 

\begin{itemize}{\rm }
    \item[{\rm(i)}] for $\alpha\in (0, \alpha^*)$, 
there exists $D^*=D^*(\alpha)>0$ such that if $D=D^*$ and $d\not=D^*$, {then} the semi-trivial steady state $(0,v^*)$ of \eqref{eq:6} is globally stable;

\item[{\rm(ii)}]  for $\alpha\ge \alpha^*$, if $D>d>0$, 
then 
$(0,v^*)$ of \eqref{eq:6} is globally stable.
\end{itemize}
\end{conjecture}
We observe that if $r$ is non-constant and $\alpha=0$,  the slower diffusion 
rate is always favored. This suggests that $D^*(\alpha)\to 0$ as $\alpha\to 0$.
Furthermore, we expect 
$D^*(\alpha)\to +\infty$ as $\alpha\to \alpha^*-$.
We refer to \cite{LouZhaoZhou, ZhaoZhou2016} for some recent developments for the case of decreasing $r$.
The general case of $r(x)$ is  studied in  \cite{LamLouLut}, where 
{sufficient conditions for the existence of evolutionarily stable strategies, or that of branching points, are found. We refer the readers to \cite{odo} for the definition of a branching point, and to Subsection \ref{S4.2} for  further discussions.}  

For the  Speirs and Gurney model \eqref{eq:3}, it was conjectured  that  there is an intermediate diffusion rate  which is optimal for  the evolution of dispersal.  This seems to be in contrast to Theorem \ref{law:2}, where the faster diffusion  is favored. Such discrepancy  is due to the different boundary conditions at the downstream end $x =L$; see \cite{LouLutscher2014,LouZhou2015, Vasilyeva2011} for further  results on the effect of boundary conditions.
{Next, we impose  the zero Dirichlet boundary 
conditions at the downstream end 
and consider the following competition model:}
 \begin{equation}\label{eq:6-a}
\left\{
\aligned
&\partial_tu=d \partial_{xx}u -\alpha \partial_xu + u (r-u-v) && (x,t)\in (0,L)\times (0, \infty),
\\
 &\partial_tv = D\partial_{xx}v -\alpha \partial_xv + v (r-u-v) && (x,t)\in (0,L)\times (0, \infty),
\\
& d \partial_xu-\alpha u= D \partial_xv-\alpha v=0 \quad &&(x,t)\in\{0\}\times (0, \infty), \\
&{u= v=0}  \quad &&(x,t)\in\{L\}\times (0, \infty), \\
& u(x,0)=u_0(x), \ v(x, 0)=v_0(x)
&& x\in(0, L).
\endaligned
\right.
\end{equation}
 As shown in Subsection \ref{S2.2}, neither  large nor small diffusion rate is beneficial to each species  in this case. A natural conjecture is: 
\begin{conjecture}\label{conjecture2} 
 There exists some $D^*(\alpha)>0$ such that if $D=D^*$ and $d\not=D^*$, the semi-trivial steady state $(0,v^*)$ of \eqref{eq:6-a},
 whenever it exists, is globally stable.
\end{conjecture}

For the case of  weak advection ($\alpha$  small), 
we conjecture that the approaches  developed in  \cite{LamLou2014-2, LamLou2014-1} may be  applicable 
 to yield the existence of $D^*$ for small $\alpha$, satisfying  $D^*(\alpha)\to 0$
 as $\alpha\to 0$. 
 However,  the proof  remains quite challenging  for general $\alpha$.
 Recently it was shown in  \cite{XiangFang}  that such diffusion rate $D^*$  exists and is unique
 for two-patch models: The semi-trivial steady state $(0, v^*)$ 
 is locally stable when $D=D^*$ and $d\not=D^*$. However, the global stability of $(0, v^*)$ remains unclear. 
 {For multi-patch models, it remains open to determine the existence of $D^*$.}

 

\section{\bf Competition models with directed movement}\label{S4}

Belgacem and Cosner \cite{BC1995} 
argued that in the spatially heterogeneous environment, the species may have a tendency to move upward along the resource gradient, in addition to random
diffusion. In this connection, 
they added an advection term to \eqref{eq:1} and  considered the reaction-diffusion-advection model
 \begin{equation}\label{eq:advection-1}
\left\{
\aligned
&\partial_t u=\nabla\cdot[d\nabla u-\alpha u\nabla m]+ u(m-u) \quad &&(x,t)\in\Omega\times (0, \infty),
\\
&\frac{\partial u}{\partial\nu}-\alpha u\frac{\partial m}{\partial\nu}=0\quad &&(x,t)\in\partial\Omega\times (0, \infty),
\\
&u(x,0)=u_0(x) \quad &&x\in\Omega.
\endaligned
\right.
\end{equation}
Here {$\alpha \geq 0$} is a parameter  which  measures the tendency of the population to move upward along the gradient of the resource $m(x)$. Belgacem and Cosner  \cite{BC1995} and subsequent work \cite{CL2003} indicate that for a single species, the advection along the resource gradient can be beneficial to the persistence of the single species in  convex habitats, but not necessary for non-convex $\Omega$. In particular, when $m$ is positive somewhere,  the species persists whenever the advection rate $\alpha$ is large. The biological explanation is that the species can utilize the large advection to concentrate at places with locally most abundant  resources. 

It is natural to enquire whether the advection along the resource gradient can confer competition advantage. 
To this end,  Cantrell et al. \cite{CCL2006} proposed the following model:
 \begin{equation}\label{eq:advection-2}
\left\{
\aligned
& \partial_t u=\nabla\cdot[d\nabla u-\alpha u\nabla m]+ u(m-u-v)\quad &&(x,t)\in\Omega\times (0, \infty),
\\
&\partial_t v=D\Delta v+ v(m-u-v) \quad &&(x,t)\in\Omega\times (0, \infty),
\\
&d\frac{\partial u}{\partial\nu}-\alpha u\frac{\partial m}{\partial\nu}=\frac{\partial v}{\partial\nu}=0\quad &&(x,t)\in\partial\Omega\times (0, \infty),
\\
&u(x,0)=u_0(x), \ v(x, 0)=v_0(x) \quad &&x\in\Omega.
\endaligned
\right.
\end{equation}

When $m$ is non-constant and $\int_\Omega m\geq 0$, \eqref{eq:advection-2}  has two semi-trivial states, denoted
by $(u^*,0)$ and $(0, v^*)$, for every $d, D>0$ and $\alpha\geq 0$ (see \cite{CL2003}),
where  $u^*$ is the unique positive steady state of \eqref{eq:advection-1}. It is shown in \cite{CCL2007} that if $d=D$ and $\Omega$ is convex, then the semi-trivial state $(u^*, 0)$ is globally stable for small $\alpha>0$. This implies that for two species which differ only in their dispersal mechanisms, the weak advection can confer some competition advantages, provided that the underlying habitat is convex.  However,  for some  non-convex  $\Omega$, it is also proved in \cite{CCL2007} that  the opposite  case may happen for  proper  $m (x)$: the weak advection can actually  put the species at a disadvantage.

Does the species have a competitive advantage  if the advection is strong? 
The following result of \cite{ALL2015} might seem surprising at the first look. It was first proved in \cite{CCL2007} under the additional assumption that the set of critical points is of measure zero and has at least one strict local maximum point.

\begin{theorem}\label{adevtion-coexistence}
Suppose that $m\in C^2({\overline\Omega})$ is {positive} and non-constant. Then for $\frac{\alpha}{d} \geq \frac{1}{\min_{\overline\Omega} m}$,
both semi-trivial steady states $(u^*, 0)$ and $(0,v^*)$ are unstable, and \eqref{eq:advection-2} has at least one stable coexistence state. 
\end{theorem}
Theorem \ref{adevtion-coexistence} implies that the two competing species are able to coexist when the advection is strong, which  is in stark contrast with the case of weak advection. 
What is preventing the species $u$ to exclude species $v$ by strong advection?
A series of work suggests that as $\alpha$ increases,  {species $u$ resides only 
at the set $\mathfrak{M}$ of local maximum points of $m(x)$, and retreats from $\Omega\setminus\mathfrak{M}$.
If $m$ is positive, then that leaves ample room for species $v$ to survive in $\Omega \setminus \mathfrak{M}$, i.e. places with less resources \cite{CCL2007, CL2008, Lam2011,  Lam2012, LamNi2010}. This is proved in Theorem \ref{adevtion-coexistence}. 
When $m(x)$ changes sign, however, the overall leftover resources in $\Omega\setminus\mathfrak{M}$ might be poor and $u$ may competitively exclude $v$ by strong advection.
\begin{theorem}[\cite{LamNi2014}]\label{thmlamni}
Suppose that $m\in C^2({\overline\Omega})$ and it is {sign-changing}. 
\begin{itemize}
    \item[{\rm(a)}] If $\int_{\Omega \setminus \mathfrak{M}} m >0$, then for any $d,D>0$, $u$ and $v$ coexist for large $\alpha$.
    \item[{\rm(b)}] If $\int_{\Omega \setminus \mathfrak{M}} m <0$, then for some $d,D>0$, $u$ competitively excludes $v$ for large $\alpha$.
\end{itemize}
\end{theorem}}
Theorem \ref{adevtion-coexistence} also reveals  that the two competitors modeled by  \eqref{eq:advection-2} can achieve coexistence by different dispersal strategies including random diffusion and directed  advection towards  resources,
which suggests a new mechanism for coexistence of two or more competing species. 
In \cite{LM2012} we apply the similar ideas to discuss the possibility of coexistence for three-species competition systems with the random diffusion and directed  advection. The analysis of three-species models turns out to be much more difficult. 
An obvious challenge is that systems of three competing species, unlike two-species competition model, are not order-preserving \cite{Smith}. 

As mentioned  earlier,  if the advection rate $\alpha$ is large, then the species $u$  will concentrate near some local maximum points of $m$.
Further studies  \cite{Lam2012, LamNi2010} showed that the concentration  positions are precisely  
those local maximum points of $m$
where the function $m-v^*$ is also strictly positive, {by establishing the uniform $L^\infty$ estimate and a new Liouville-type result}. We note that $m$ is the common resource for both species $u$ and $v$, and $m-v^*$ is the  resource available for species $u$ when the species  $v$  reaches its equilibrium. It is observed in \cite[Rmk. 1.4]{LamNi2010} that for some positive local maximum points of $m$, function  $m-v^*$ can be negative, i.e. some apparently favorable places are actually dangerous ecological traps for the species $u$.
This suggests  that both  distribution of  resources and competitors should be taken into account when studying the directed movement of species along the resource gradient. 

\subsection{Global bifurcation diagrams}
\begin{figure}
  \centering
  \includegraphics[width=0.5\linewidth]{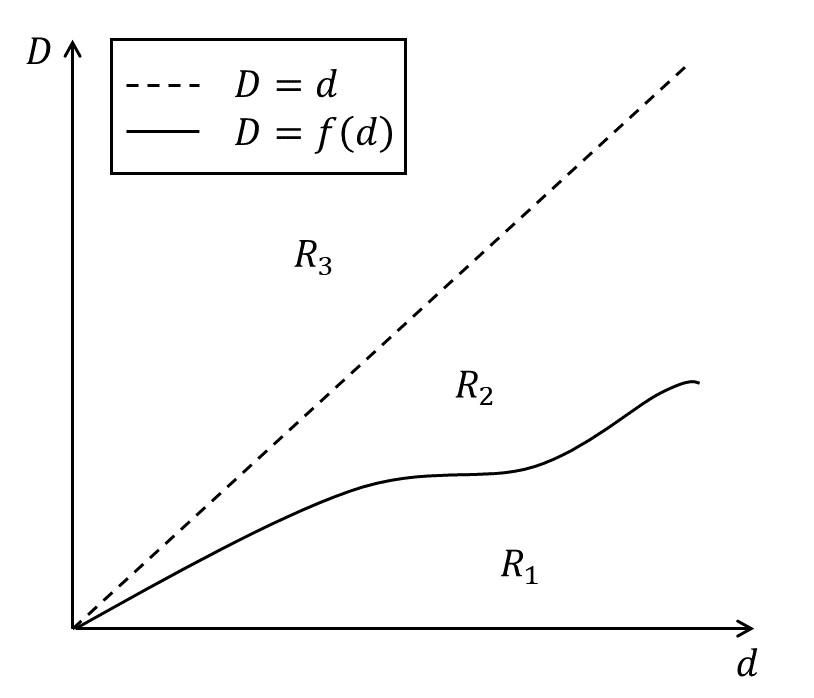}
  \caption{Decomposition of the $d$-$D$ parameter space into three regions: $R_1$ and $R_2$ are separated by the curve $D = f(d)$. The regions $R_2$ and $R_3$ are separated by the line $D=d$. The bifurcation diagrams for each of these regions, with $\alpha$ as the bifurcation parameter, are given in Figure \ref{mams1}. (First published in Mem. Amer. Math. Soc. 245 (Spring 2017), published by the American Mathematical Society. @ 2017 American Mathematical Society.)}\label{mams2}
 \end{figure}

\begin{figure}
  \centering
  \includegraphics[width=0.5\linewidth]{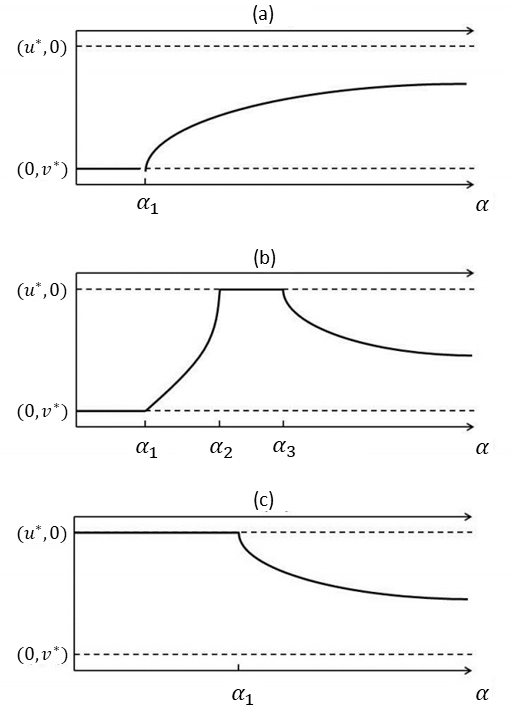}
  \caption{Bifurcation results from \cite{ALL2015}, where $d,D$ are fixed and $\alpha$ is used as the bifurcation parameter. 
  The three figures, from top to bottom, are the bifurcation diagrams concerning the parameter regions $R_1, R_2$, and $R_3$ (respectively, $D < f(d)$, $f(d) < D < d$ and $D > d$) in Figure \ref{mams2}.  (First published in Mem. Amer. Math. Soc. 245 (Spring 2017), published by the American Mathematical Society. © 2017 American Mathematical Society.) }\label{mams1}
 \end{figure}

The above results show that when  the advection is weak, the two species cannot coexist, and  the advection confers some
competitive advantage for species $u$; when  the advection is strong, the two species coexist for any diffusion rates. 
What happens for general advection rates? For instance,  assuming $d>D$, when we increase the advection from weak to strong, how does {that affect} the dynamics of \eqref{eq:advection-2}? If $\alpha=0$, the result of Dockery et al. \cite{DHMP} says that the semi-trivial state $(0, v^*)$ is globally stable; if  $\alpha$ is sufficiently large, Theorem \ref{adevtion-coexistence} shows that both $(u^*, 0)$ and $(0,v^*)$ are unstable and the coexistence happens.  A possible conjecture is there exists some critical advection rate $\alpha^*>0$ such that for $\alpha<\alpha^*$, $(0,v^*)$ is globally stable, while for $\alpha>\alpha^*$, \eqref{eq:advection-2} has a unique positive steady state which is globally stable. However this conjecture fails in general and a more complicated bifurcation diagram of positive steady states occurs. To be more precise, Averill et al. \cite{ALL2015} showed that  under some
conditions on $m(x)$, 
there exists a function $f(x): [0,\infty)\to [0,\infty)$ satisfying $f(0)=0$ and $0<f(x)<x$ such that the $d$-$D$ parameter space is divided into three regions $R_1, R_2, R_3$, as shown in Figure \ref{mams2}. 
Using $\alpha$ as a bifurcation parameter, we see that the following three situations occur:
\begin{enumerate}
\item for $D<f(d)$, there exists some critical value $\alpha_1>0$ such that as $\alpha<\alpha_1$,  the semi-trivial state $(0,v^*)$ is stable but $(u^*,0)$ is unstable, while as $\alpha>\alpha_1$, both $(u^*,0)$ and  $(0,v^*)$ are unstable so that \eqref{eq:advection-2} has at least one positive steady state.  By the bifurcation theory,  there is a branch of  positive steady states, denoted by $\Gamma_1$, which bifurcates from $(0, v^*)$ at $\alpha=\alpha_1$ and can be extended to infinity in $\alpha$. This is illustrated in Figure \ref{mams1}(a).
It is unclear whether $(0, v^*) $ is globally stable when $\alpha <\alpha_1$, and whether the positive steady state  is unique and  globally stable when $\alpha>\alpha_1$;
\smallskip

\item for $f(d)<D<d$, there exist three critical values  $0<\alpha_1<\alpha_2<\alpha_3$ such that {\rm(i)} as $\alpha<\alpha_1$,  $(0,v^*)$ is stable but $(u^*,0)$ is unstable;  {\rm(ii)} as $\alpha_2<\alpha<\alpha_3$, $(u^*, 0)$ is stable but $(0,v^*)$ is unstable;  {\rm(iii)} as $\alpha\in (\alpha_1,\alpha_2)\cup (\alpha_3, \infty)$, both $(u^*, 0)$ and $(0,v^*)$ are unstable so that  \eqref{eq:advection-2} has at least one positive steady state.  Similarly, by the bifurcation theory, there is still a branch of   positive steady states  bifurcating from $(0, v^*)$ at $\alpha= \alpha_1 $, which is  denoted as $\Gamma_2$.  The difference is that the branch  $\Gamma_2$ cannot be extended to $\alpha=\infty$, but is connected to $ (u^*,0) $ at $\alpha =\alpha_2$.  Interestingly, when $\alpha\in(\alpha_2, \alpha_3) $,  the semi-trivial state $(u^*, 0)$ is stable, which indicates that the advection can indeed be beneficial to the species in this case!  When  $\alpha>\alpha_3$, the branch of a  positive steady state,  denoted as $\Gamma_3$,  bifurcates from $ (u^*, 0)$ at $\alpha=\alpha_3 $, and can be extended to $\alpha=\infty$. We conjecture that as $\alpha<\alpha_1$, $(0,v^*)$ is globally stable; as $\alpha_2<\alpha<\alpha_3$, 
$(u^*, 0)$ is globally stable; as $\alpha\in (\alpha_1,\alpha_2)\cup (\alpha_3, \infty)$, 
\eqref{eq:advection-2} admits a unique positive steady state which is globally stable. This is illustrated in Figure \ref{mams1}(b).

\smallskip 

\item for $D>d$, there exists some critical value $\alpha_1>0$ such that as $\alpha<\alpha_1$, $(u^*, 0)$ is stable, while as $\alpha>\alpha_1$, both $(u^*, 0)$ and $(0,v^*)$ are unstable so that  \eqref{eq:advection-2} has at least one positive steady state.
Similarly,  a branch   of  positive steady states bifurcates from $(u^*, 0)$ at $\alpha=\alpha_1$, denoted as $\Gamma_4$, which can be extended to $\alpha=\infty$. 
It is conjectured that $(u^*, 0)$ is globally stable when $\alpha<\alpha_1$, and the positive steady state is unique when $\alpha>\alpha_1$. This is illustrated in Figure \ref{mams1}(c).
\end{enumerate}

How, then, is the branch $\Gamma_1$ 
transformed to $\Gamma_2$ and $\Gamma_3$, and finally to $\Gamma_4$,
as the diffusion rate $D$ varies,
with $d$ being fixed? 
The result in \cite{ALL2015} shows that in the critical situation $D=f(d)$, the branch $\Gamma_1$ is connected to 
the semi-trivial state $ (u^*, 0) $ at  $ \alpha=\alpha^*$ for some $\alpha^*>0$. Once $D> f(d)$,  $\alpha^*$ is  split into two different points $\alpha_2$ and $\alpha_3 $, so that $\Gamma_1$ is split into two disjoint branches  $\Gamma_2$ and $\Gamma_3$.
As $D$ approaches $d$ from below, the bounded branch $\Gamma_2$ will approach $\alpha=0$ and vanishes; meanwhile  $\Gamma_3$ approaches $\Gamma_4$.

By these discussions, we explained that for the models with multiple parameters, it is possible to  apply the bifurcation theory to describe the structure of  positive steady states.  Of course, it remains challenging  to completely characterize the global dynamics of \eqref{eq:advection-2}, and many interesting mathematical questions, e.g. uniqueness and global stability, are left open \cite{ALL2015}. 
One particular conjecture is that, for $\alpha$ sufficiently large, the positive steady state of the competition system \eqref{eq:advection-2} is unique and thus globally asymptotically stable. One way of proving it is to determine the precise limiting profile and show that each positive steady state is linearly stable, whenever it exists. One can then conclude the uniqueness and global asymptotic stability by the theory of monotone dynamical system \cite{Hess}. This was carried out for the case when all critical points of $m$ are non-degenerate, and that the value of $m$ is constant on the set of all local maximum points \cite[Thm. 8.1]{lamthesis}. The general case remains open.

It is  illustrated here   that  proper advection can confer some competition advantages. In  next subsection,  we will use evolutionary game theory to further show that proper advection can be an  evolutionary stable strategy.

 

\subsection{\bf Evolution of biased movement}\label{S4.2}

In this subsection, we adopt the viewpoint of adaptive dynamics \cite{odo,Geritz1998} to consider the evolution of dispersal and the ecological impact of dispersal in spatially heterogeneous environments. An important concept in adaptive dynamics is  Evolutionarily Stable Strategy (ESS), known also as Evolutionarily Steady Strategy, which was introduced by 
Maynard Smith and Price in the seminal paper \cite{ESS}. A strategy is said to be evolutionarily stable if a population using it cannot be invaded by any small population using a different strategy. ESS is a very general concept. 
In this subsection we aim to connect it  with the local stability of the semi-trivial steady 
states discussed earlier.

An interesting question is  what kind of advection rate will be evolutionarily stable, by regarding the rates of advection along the environmental gradient as movement strategies of the species. For this purpose, we consider the following model:
\begin{equation}\label{eq:advection-3}
\left\{
\aligned
& \partial_t u=\nabla\cdot[d\nabla u-\alpha u\nabla m]+ u(m-u-v)\quad &&(x,t)\in\Omega\times (0, \infty),
\\
&\partial_t v=\nabla\cdot[d\nabla v-\beta v\nabla m]+ v(m-u-v) \quad &&(x,t)\in\Omega\times (0, \infty),
\\
&d\frac{\partial u}{\partial\nu}-\alpha u\frac{\partial m}{\partial\nu}=
d\frac{\partial v}{\partial\nu}-\beta v\frac{\partial m}{\partial\nu}=0\quad &&(x,t)\in\partial\Omega\times (0, \infty),
\\
&u(x,0)=u_0(x), \ v(x, 0)=v_0(x) \quad &&x\in\Omega,
\endaligned
\right.
\end{equation}
where $u$ and $v$ are the densities of resident and mutant populations, and the two species differ only in their advection rates, given by the non-negative parameters $\alpha$ and $\beta$, respectively. 

The mutant population $v$, which is assumed to be small initially, can invade successfully if and only if the semi-trivial steady state $(u^*, 0)$ is unstable. 
Here $u^*$ denotes the unique positive solution of 
\begin{equation}\label{eq:advection-4}
\left\{
\aligned
& \nabla\cdot[d\nabla u-\alpha u\nabla m]+ u(m-u)=0\quad &&x\in\Omega,
\\
&d\frac{\partial u}{\partial\nu}-\alpha u\frac{\partial m}{\partial\nu}=0
\quad &&x\in\partial\Omega.
\endaligned
\right.
\end{equation}

A natural question is whether there is a critical value  $\alpha^*>0$ such that $(u^*, 0)$ is a stable steady state of \eqref{eq:advection-3} whenever $\alpha=\alpha^*$ and $\beta\not=\alpha^*$. 
{We call such $\alpha^*$, if it exists, an ESS}. To find such a magical strategy, we note that  the stability of  $(u^*, 0)$ is determined by the principal eigenvalue,
denoted by $\Lambda_1$, of the problem 
\begin{equation}\label{eq:advection-5}
\left\{
\aligned
& \nabla\cdot[d\nabla \varphi-\beta \varphi\nabla m]+ \varphi(m-u^*)+\lambda\varphi=0\quad &&x\in\Omega,
\\
&d\frac{\partial \varphi}{\partial\nu}-\beta \varphi\frac{\partial m}{\partial\nu}=0
\quad &&x\in\partial\Omega.
\endaligned
\right.
\end{equation}
It is classical that the set of eigenvalues of \eqref{eq:advection-5} is real and bounded from below. The principal eigenvalue $\Lambda_1$ refers to the smallest eigenvalue. It is simple and corresponds to a strictly positive eigenfunction in $\overline\Omega$; see \cite{CC2003,CL2008}. Since $u^*$ is an implicit 
function of  $\alpha$, we regard $\Lambda_1 = \Lambda_1(\alpha,\beta)$. 
Observe that  $\Lambda_1=0$ for $\alpha=\beta$, since the principal eigenfunction $\varphi$ is simply a constant multiple of $u^*$ in such a case. 

 It is well known that if $\Lambda_1>0$, then $(u^*, 0)$ is stable; otherwise,  $(u^*, 0)$ is unstable. Hence if  $\alpha^*$ is an ESS, then $\Lambda_1(\alpha^*, \beta)\ge 0$ for all $\beta$, which together with $\Lambda_1(\alpha^*, \alpha^*)= 0$ implies that
 \begin{equation}\label{eq:advection-6}
\frac{\partial\Lambda_1}{\partial\beta}\Big|_{\alpha=\beta=\alpha^*}=0.
\end{equation}
We call any strategy $\alpha^*$  satisfying \eqref{eq:advection-6} an evolutionarily singular strategy, 
hence any ESS is automatically an evolutionarily singular strategy.  For  the existence of evolutionarily singular strategies,  we have the following result:
\begin{theorem}[{\rm \cite[Thm. 2.2]{LamLou2014-2}}] \label{thm:singularstrategy}
Suppose that $m\in C^2({\overline\Omega})$ is strictly positive in ${\overline\Omega}$ and 
$\max_{{\overline\Omega}}m/\min_{{\overline\Omega}} m\le 3+2\sqrt{2}$. Given any 
 $A>1/\min_{{\overline\Omega}} m$, for all small
positive $d$, there is exactly one evolutionarily singular strategy  $\alpha^*\in (0, dA)$. 
\end{theorem} 

The condition $\max_{{\overline\Omega}}m/\min_{{\overline\Omega}} m\le 3+2\sqrt{2}$ is sharp, and Theorem \ref{thm:singularstrategy} may fail for general $m$ \cite[Thm. 2.6]{LamLou2014-2}. 
{When $\Omega$ is one-dimensional, the evolutionarily singular strategy found in Theorem \ref{thm:singularstrategy} is also unique in the  whole interval $[0, \infty)$, provided $d>0$ is small enough \cite[Cor. 2.3]{LamLou2014-2}. We conjecture that the uniqueness holds for $\Omega$ of higher dimensions as well.} For intermediate or large $d$, however, the uniqueness of singular strategies remains open; see \cite{CCL2019} for some recent progress.
We conjecture that for any  evolutionarily singular strategy $\alpha^*=\alpha^*(d)$, it holds that, for all $d$,  
\begin{equation}\label{eq:singular-11}
\min_{{\overline\Omega}} m<\frac{d}{\alpha^*}<\max_{{\overline\Omega}} m.
\end{equation}

If this estimate were true, it suggests that a balanced combination of diffusion and advection might be a better movement strategy for populations.  The following result shows that this is indeed the case:
 \begin{theorem}[{\rm\cite[Thm. 2.5]{LamLou2014-2}}]
 Suppose that $\Omega$ is convex and $m\in C^2({\overline\Omega})$ is strictly positive in ${\overline\Omega}$. Assume that  $\|\nabla (\ln m)\|_{L^\infty(\Omega)}\le \alpha_0/{\rm diam}(\Omega)$, where 
$\alpha_0\approx 0.8$ and ${\rm diam}(\Omega)$ is the diameter of $\Omega$. 
Then for small $d$, the strategy $\alpha^*$ given by  Theorem {\rm\ref{thm:singularstrategy}} is a local ESS, i.e. 
there exists some small $\delta>0$ such that 
$\Lambda_1(\alpha^*, \beta)>0$ for all $\beta\in(\alpha^*-\delta, \alpha)\cup (\alpha, \alpha^*+\delta)$.
 \end{theorem}\label{thm:ess-2}

Some other work  also considered model \eqref{eq:advection-3} and  related models. For instance, in \cite{LamLou2014-1} we  discussed the evolution of diffusion rate in an advective environment. We refer to \cite{CCL2019, CHL2008, CLL2012, Cosner2013, Cosner2014,  Gejji, HL2009} for further discussions.

{Can we decide whether the local ESS $\alpha^*$ given in Theorem {\rm\ref{thm:singularstrategy}} is actually a global ESS? Or more generally, what can we say about the global structure of nodal set for the function  $\Lambda_1=\Lambda_1(\alpha, \beta)$?} 
It is clear that the nodal set of $\Lambda_1(\alpha,\beta)$ must include the diagonal line $\alpha=\beta$, but  otherwise 
little is known. 
For the reaction-diffusion models discussed in this paper, {determining the complete structure}
of the nodal set for $\Lambda_1(\alpha, \beta)$
poses  a challenging problem, both analytically and numerically \cite{GHLL}. To this end,  some asymptotic behaviors  of $\Lambda_1$ were studied in \cite{CL2008, CL2012, LL2019, PZ2017} in order to have a better understanding for this problem.  Recently, the  asymptotic behaviors for  the principal eigenvalues of time-periodic parabolic operators were also studied,
and we refer to \cite{Hess, Hutson2000,Hutson2001, LLPZ20191, LLPZ20201,LLPZ20202,Peng2015}. 
{These new developments of the related linear PDE theory may open up a new venue to study the evolution of dispersal in
spatio-temporally varying environments. }

The evolutionary dynamics, e.g. existence of ESS, have various consequences on the two-species competition dynamics. This connection has been investigated in \cite{CCL2017} for a broad class of model including reaction-diffusion equations and nonlocal diffusion equations. It was shown that frequently the ESS playing species can competitively exclude the species playing a different but nearby strategy, regardless of the initial condition. In the jargon of adaptive dynamics, this result says that ESS implies Neighborhood Invader Strategy. 
This is quite unexpected, since an ESS by its very definition only concerns the local stability of the semi-trivial steady state. The main tool in \cite{CCL2017} is Hadamard's graph transform method.

\section{\bf Continuous trait models}\label{S5}

Regarding the advection rates as  strategies of the species, in Subsection \ref{S4.2} we {investigated} the course of evolution for two-species competition model \eqref{eq:advection-3}. Therein the strategies $\{\alpha,\beta\}$ form a discrete set, and  the competition and mutation occur independently of each other. In this section, we consider population models with continuous traits, in which the competition and mutation can occur simultaneously.

\subsection{\bf Mutation-selection model}\label{S4.2.1}

To study the selection of random
dispersal rate {in} multi-species competition, Dockery et al. \cite{DHMP} considered the following model:
\begin{equation}\label{liu01}
    \partial_t u_i=\xi_i\Delta u_i+ u_i\left[m(x)-\sum_{j=1}^K u_j\right] +\sum_{j=1}^K M_{ij} u_j, 
\end{equation}
where $K$ interacting species, with densities $u_i(x,t)$, $1\leq i\leq K$, compete for a common resource $m(x)$ and differ only in their diffusion rates $\xi_i$. These species are subject to mutation with switching rate $M_{ij}$ from phenotype $j$ to $i$. 
For $K=2$ and $M_{ij}=0$, a well-known result is that the slower diffusing species always win in the competition, but it remains open to determine the global dynamics for the case ~$K\geq 3$.

Next, we extend the discrete set of strategies $\{\xi_i\}_{i=1}^K$ in \eqref{liu01} to a continuum set $(\underline{\xi},\overline{\xi})$. By viewing the population $u(x,\xi,t)$ as being structured by space, trait and time (i.e. $u(\cdot,\xi,t)$ gives the spatial distribution of the individuals playing strategy $\xi$ at time $t$), 
we obtain the following mutation-selection model:
\begin{equation}\label{eq:8}
\left\{
\aligned
&\partial_t u=\xi \Delta_x u
+ u \left[m(x) -  \hat{u}(x,t)
\right]+\epsilon^2 \partial_{\xi\xi} u
&& \hbox{ for }\,x\in\Omega, \xi\in (\underline\xi, \overline\xi), t>0,\\
&{\nu(x) \cdot \nabla_x u =0}
&&\hbox{ for }\,x\in\partial \Omega,  \xi\in (\underline\xi, \overline\xi), t>0, \\
&\partial_{\xi} u=0 && \hbox{ for }\,x\in\Omega, \xi\in \{\underline\xi, \overline\xi\}, t>0,
\\
&u(x,\xi, 0)=u_0(x,\xi)
&& \hbox{ for }\,x\in\Omega, \xi\in (\underline\xi, \overline\xi).
\endaligned
\right.
\end{equation}
Here the parameters $\underline\xi, \overline\xi$, which are assumed to be positive, determine the range of continuous strategy $\xi$, 
and
$$ \hat{u}(x,t):=\int_{\underline\xi}^{\overline\xi} u(x,\xi', t)\,\mathrm{d}\xi'$$ represents the total  population density at location $x$ and time $t$, and the term $\epsilon^2 \partial_{\xi\xi} u$ accounts for
the rare mutation, which is modeled by a diffusion with covariance $\sqrt{2}\epsilon$. We refer to \cite{CF2006} for a derivation of this type of equations from individual based {on} stochastic models. 
Model \eqref{eq:8} is an integro-PDE model describing the diffusion, competition and mutation of a  population with continuous traits. 

If $\int_{\Omega}m(x)\mathrm{d}x>0$, then it is shown in  \cite{Lam2017} that \eqref{eq:8} admits a unique positive steady state $u_\epsilon(x, \xi)$ for small $\epsilon$, which is locally asymptotically stable and satisfies 
\begin{equation}\label{eq:8-1}
\left\{
\aligned
&\xi \Delta_x u + \epsilon^2 \partial_{\xi\xi} u  + u \left[m(x)
- \hat{u}(x)
\right] = 0 
\ && (x,\xi)\in\Omega \times (\underline\xi, \overline\xi),
\\
&{\nu(x) \cdot \nabla_x u=0}
 \ &&(x,\xi)\in\partial \Omega \times (\underline\xi, \overline\xi), 
\\
&\partial_{\xi} u=0 
\ && (x,\xi)\in\Omega \times \{\underline\xi, \overline\xi\},
\endaligned
\right.
\end{equation}
where $\hat{u}(x)=\int_{\underline\xi}^{\overline\xi} u(x,\xi')\mathrm{d}\xi'$. 
Observe that when $m$ is constant, i.e. $m\equiv m_0$ for some positive constant $m_0$,  we have $u_\epsilon\equiv m_0/(\overline{\xi}-\underline{\xi})$ which is globally stable for any $\epsilon>0$ \cite{Lam2017}. This suggests that {selection does not favor any particular strategy $\xi \in [\underline\xi,\overline\xi]$, and} 
the mutation $\epsilon^2 \partial_{\xi\xi} u$ has no effect on the population dynamics in the homogeneous environment. What happens in the spatially heterogeneous environment? The following result gives the answer:

\begin{theorem}[{\rm\cite{LL2017}}]\label{law:3}
Suppose that $\int_{\Omega}m(x)\mathrm{d}x>0$ and $m$ is non-constant. Then as  $\epsilon \to 0$, 
\begin{equation*}
\left\|
\epsilon^{2/3}u_\epsilon(x,\xi) -u^*(x,{\underline\xi})\eta^*\left( \frac{\xi - \underline\xi}{\epsilon^{2/3}}\right)\right\|_{L^\infty(\Omega\times (\underline{\xi}, \overline{\xi}))} \to 0, 
\end{equation*}
where $\eta^*(s)$ satisfying the following  ordinary differential equation:
\begin{equation*}
\left\{
\begin{array}{ll}
\medskip
\eta'' + (a_0 -  a_1 s)\eta = 0\quad \quad s>0,\\
\eta'(0) =\eta(+\infty)=0, \quad \, \int_0^\infty \eta(s)\,\mathrm{d}s = 1,
\end{array}
\right.
\end{equation*}
 with some positive constants $a_0$ and $ a_1$. In particular, we have as $\epsilon\to 0$, 
$$u_\epsilon(x, \xi)\to u^*(x, \underline{\xi}) \cdot \delta(\xi-\underline{\xi}) \quad\text{in distribution},$$
where $ u^*(x, \underline{\xi})$ denotes the unique positive solution of \eqref{eq:2} with $d=\underline{\xi}$, and $\delta(\cdot)$ is a Dirac mass supported at $0$.
\end{theorem}

It can be verified that for any $\xi_0\in [\underline\xi,\overline\xi]$, the limit profile $u^*(x, \xi_0)\cdot \delta(\xi-\xi_0)$ is a solution of \eqref{eq:8-1} with $\epsilon=0$. {Theorem \ref{law:3}} 
says that only the steady state corresponding to the smallest diffusion rate is preserved by 
the perturbation of $\epsilon =0$ to $0<\epsilon\ll 1$. 
Theorem \ref{law:3} implies that  when the mutation rate 
is small, the species with the smallest diffusion rate  $\underline{\xi}$ will eventually become the dominant phenotype by driving other phenotypes to extinction. This result {connects} with  the case of two-species competition:  The smallest available diffusion rate is selected. Note from  \cite{Lam2017} that the steady state $u_\epsilon$ is the unique solution of \eqref{eq:8-1} for small  $\epsilon>0$, and it is locally asymptotically stable. 
 An unsolved  problem here is the global stability of the unique positive steady state $u_\epsilon(x,\xi)$.

A {closely} related work is due to Perthame and Souganidis \cite{PS2015}, where the following mutation-selection model was considered:
\begin{equation}\label{eq:9}
\left\{
\aligned
&\xi(z) \Delta_x u + \epsilon^2 \partial_{zz} u + u (m (x) - \hat{u}) = 0 && (x,z)\in \Omega \times [0,1],\\
&{\nu(x) \cdot \nabla_x u =0}
 && (x,z)\in\partial \Omega \times [0,1],\\
&u(x,0)=u(x,1) 
&& x\in \Omega,
\endaligned
\right.
\end{equation}
where $\Omega$ is assumed to be convex and $\xi(z)$ is assumed to be positive and periodic with unit period.  For the solution $u_{\epsilon}(x, z)$ of \eqref{eq:9}, it is proved in 
\cite{PS2015} that  $u_{\epsilon}(x, z) \to u^*(x, \xi(z_*)) \cdot \delta(z - \xi(z_*))$ in distribution as $\epsilon\to 0$,
where $\xi(z_*)=\min_{0\le z\le 1} \xi(z)$. This implies that  for rare mutations, the population concentrates on a single trait associated to the smallest diffusion rate, in agreement with the result  for \eqref{eq:8-1}. See also \cite{BM2005, JS2020, MB2015} for some progress in this direction. 

\subsection{\bf Mutation-selection model with drift}\label{S4.2.2}
Similar to \eqref{eq:8}, we {introduce} the following integro-PDE model in a stream:
\begin{equation}\label{eq:10}
\left\{
\aligned
&\partial_tu=\xi \partial_{xx}u -\alpha\partial_xu + u\left[r(x) - \hat u(x,t)\right] + \epsilon^2 \partial_{\xi\xi}u&&
  \hbox{ for }\,x\in(0,L), \xi\in (\underline\xi, \overline\xi), t>0,\\
&(\xi \partial_{x}u - \alpha u)(0,\xi,t)=(\xi \partial_{x}u - \alpha u)(L,\xi,t)=0 &&  \hbox{ for }\, \xi\in  (\underline\xi, \overline\xi), t>0,\\
&u(x,\overline{\xi},t)=u(x,\underline{\xi},t)=0  &&  \hbox{ for }\,x\in(0,L)\, t>0,
\\
&u(x, \xi, 0)=u_0(x, \xi)
&&\hbox{ for }\,x\in(0,L), \xi\in (\underline\xi, \overline\xi),
\endaligned\right.
\end{equation}
where $\alpha>0$ is the advection rate and $\hat{u}(x,t):=\int_{\underline\xi}^{\overline\xi} u(x,\xi', t)\,\mathrm{d}\xi'$ denotes the total  population density at location $x$ and time $t$.

 When $r$ is a positive constant, from Theorem \ref{law:2} we conjecture that {\rm (i)} the integro-PDE \eqref{eq:10} admits at most one positive steady state, denoted by $u_\epsilon(x,\xi)$, which is  globally stable whenever it exists; {\rm (ii)} as $\epsilon\to 0$, $u_\epsilon(x, \xi)\to u^*(x, {\overline\xi}) \cdot \delta(\xi-\overline{\xi})$ in distribution sense, where $u^*(x, \overline{\xi})$ is the unique positive solution of \eqref{eq:7} with $D=\overline{\xi}$.  
 Again, this predicts that the larger 
 diffusion is selected as in 
 Theorem \ref{law:2}.

\begin{figure}
  \centering
  \includegraphics[width=0.8\linewidth]{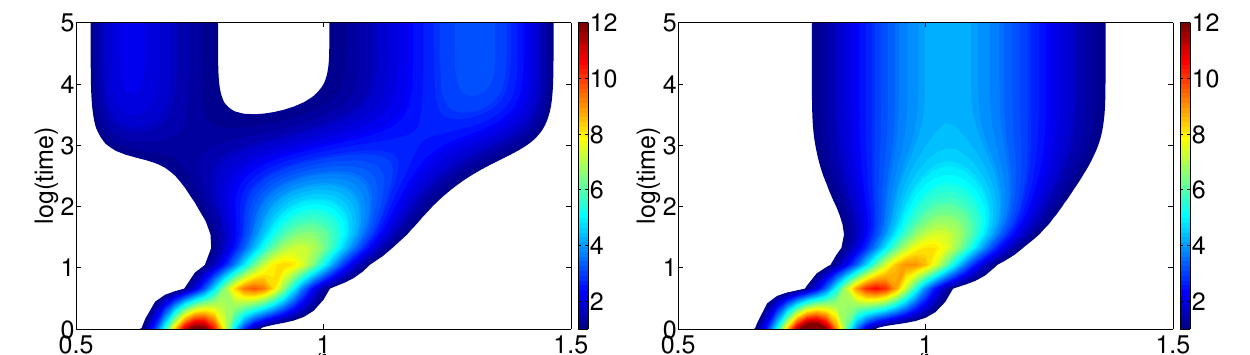}
  \caption{Contour plot of $\int u(x, \xi, t)\mathrm{d}x$
 as a function of $\xi$ and time (log(time) for vertical axis) for $a= \frac{1}{4}$
(left) and $a=-\frac{1}{4}$
(right). 
{(First published in [W. Hao et al., Indiana Univ. Math. J. 68 (2019), 881-923.].)}
}\label{liufig2}
 \end{figure}

 To study the dynamics of \eqref{eq:10} for non-constant $r$,  we  carried out some numerical simulations in  \cite{HLL2019}, where the corresponding parameters were selected as follows:
\begin{equation}\label{eq:t1}
L=1,\,\,\, \alpha=1, \,\,\, \underline\xi=0.5, \,\,\, \overline\xi=1.5,\,\,\,r(x)=e^{(1-a)x+ax^2}, \,\,\, \epsilon=10^{-3}.
\end{equation}
We take initial conditions in the form of one Dirac mass on the phenotypic space, and investigate their evolution
for $a=\pm \frac{1}{4}$. The numerical results are presented in Figure \ref{liufig2}.


\begin{itemize}
    \item[{\rm(i)}]  The  numerical result in  Figure \ref{liufig2} (right) shows that when $a=\frac{1}{4}$,
the steady state of \eqref{eq:10} concentrates on the trait $\xi^*\approx 1.1$ in the limit of rare mutation, which suggests that the species  adopting the strategy  $\xi^*$ persists  and other species will disappear in the competition. This suggests that $\xi^*$
is an ESS, which is in contrast with the assertion that ``faster diffusion  is more favorable" for the homogeneous environment. 

\smallskip 

\item[{\rm(ii)}]  
More interestingly, when $a=-\frac{1}{4}$,  Figure \ref{liufig2} (left) shows that  the steady state of \eqref{eq:10} concentrates on 
two different traits, so that two species adopting different strategies {form a coalition that is able to resist invasion by any of the remaining species}! This phenomenon is
called ``evolutionary branching" in
adaptive dynamics, and  
it is also {clearly} different from the homogeneous case.
\end{itemize}

To explain these phenomena mathematically, we investigate the asymptotic behaviors of steady state  $u_\epsilon(x, \xi)$ for \eqref{eq:10}
as $\epsilon\to 0$. To this end, we first consider the two-species competition model \eqref{eq:6} and denote $\lambda = \lambda(\xi_1,\xi_2)$ as the principal eigenvalue of the problem
\begin{equation*}
\left\{
\aligned
&\xi_2 \partial_{xx}\varphi -\alpha \partial_{x}\varphi + \left[r - {u^*(\cdot, \xi_1)}\right]\varphi +\lambda\varphi=0  &&  x \in (0,L),\\
&\xi_2 \partial_{x}\varphi - \alpha\varphi=0 &&  x = 0,L,
\endaligned
\right.
\end{equation*}
where {$x\mapsto u^*(x, \xi_1)$} is the unique positive solution of \eqref{eq:7} with $D=\xi_1$. 
In the adaptive dynamics framework,  $\lambda(\xi_1,\xi_2)$ {is called the} invasion fitness, and 
an invader with trait $\xi_2$ can (resp. cannot) invade an established trait $\xi_1$ at equilibrium when rare if $\lambda(\xi_1,\xi_2)<0$ (resp. $\lambda(\xi_1,\xi_2)>0$).
The limiting profiles of $u_\epsilon(x, \xi)$ depends critically on the behavior of $\lambda(\xi_1,\xi_2)$ 
near $\xi_2=\xi_1$. Note that $\lambda(\xi,\xi)\equiv 0$ for all $\xi>0$.  We first consider the  case of $\partial_{\xi_2}\lambda(\xi,\xi) \neq 0$, where the following result holds:
\begin{theorem}[{\rm\cite{HLL2019}}]\label{law:4}
Suppose that $\partial_{\xi_2}\lambda(\xi,\xi)>0$ for all $\xi \in (\underline{\xi}, \overline{\xi})$, and $\overline{\xi}-\underline{\xi}$ is sufficiently small. Then any positive steady state $u_\epsilon$ of \eqref{eq:10} satisfies as $\epsilon \to 0$, 
$$
u_\epsilon(x,\xi) \to \delta(\xi - \underline\xi)\cdot u^*(x, \underline\xi) \quad\text{in distribution}, 
$$
where $x\mapsto u^*(x,\underline\xi)$ is the unique positive solution of \eqref{eq:7} with $D=\underline{\xi}$.
\end{theorem}

Theorem \ref{law:4} implies that for the case  $\partial_{\xi_2}\lambda(\xi,\xi)>0$, the slowest diffusion will be selected in the competition. 
Similarly, we may conclude that if $\partial_{\xi_2}\lambda(\xi,\xi)<0$ for all $\xi \in (\underline{\xi}, \overline{\xi})$, then the steady state $u_\epsilon$ concentrates on the trait $\xi=\overline\xi$, so that  the fastest diffusion will be selected,  coinciding  with the observation in Theorem \ref{law:2}. These results show that if $\partial_{\xi_2}\lambda\neq 0$, it gives rise to a single Dirac-concentration at one of the two extreme traits, depending upon the sign of $\partial_{\xi_2}\lambda$.

To study the case when $\partial_{\xi_2}\lambda$ vanishes somewhere in $(\underline{\xi}, \overline{\xi})$, we introduce the convergence stable strategy $\hat{\xi}$ (see \cite{EM1981}),  which is characterized by the following relations:
\begin{equation}\label{eq:csss}
\begin{array}{l}
\partial_{\xi_2}\lambda(\hat\xi, \hat\xi) = 0 \quad\text{and}\quad\frac{\rm{d}}{\rm{d}s}\left[ \partial_{\xi_2}\lambda(s,s)\right]\Big|_{s=\hat\xi} > 0.
\end{array}
\end{equation}
This leads to two generic cases: The convergence stable strategy $\hat{\xi}$ is {an ESS} if  $\partial^2 _{\xi_2}\lambda(\hat \xi,\hat \xi)>0$, and is a  branching point (BP) if  $\partial^2 _{\xi_2}\lambda(\hat \xi,\hat \xi)<0$.
 
 In the case when the ESS exists, there is concentration in the mutation-selection model.
 \begin{theorem}{\rm\cite{HLL2019}}\label{law:5}
Suppose that $\hat\xi$ is an $\mathrm{ESS}$ {and \eqref{eq:csss} holds}. If $\hat\xi \in (\underline{\xi}, \overline{\xi})$ and $\overline{\xi}-\underline{\xi}$ is sufficiently small, then as $\epsilon \to 0$,  any positive steady state $u_\epsilon$ of \eqref{eq:10} satisfies
$$
u_\epsilon(x,\xi) \to u^*(x, \hat\xi)\cdot \delta(\xi - \hat\xi) \quad\text{in distribution}. 
$$
 \end{theorem}

Theorem \ref{law:5} indicates that {if there is an ESS in $(\underline{\xi}, \overline{\xi})$, then the phenotype adopting the ESS dominates the competition in the limit of rare mutation. This} corresponds to the numerical result in Figure \ref{liufig2} (right) with 
$\hat\xi\approx 1.1$. To further understand Theorem \ref{law:5}, we consider a special case of \eqref{eq:6} where $r(x)=be^{cx}$ for some positive constants $b$ and $c$. In this case, it is shown in \cite{ALM} that  the semi-trivial state $(u^*(x,\xi_1),0)$ is globally stable when $\xi_1=\alpha/c$ and $\xi_2\not=\alpha/c$, i.e. the  intermediate trait $\hat\xi=\alpha/c$ is an ESS. We can offer some intuitive reasoning for  this result by employing  the theory of ideal free distribution \cite{FL1970}:
When $r(x)=be^{cx}$, if $\xi_1=\alpha/c$,  it can be verified from  \eqref{eq:7} that
$u^*(x,\xi_1)\equiv be^{cx}\equiv r(x),$
so that species $u$ attains the ideal free distribution, and thus $\alpha/c$ being an ESS is  natural. As $c\to 0$, i.e. the environment tends to be homogeneous, we find 
$\hat\xi=\alpha/c$ approaches $ +\infty$, which is in agreement with the prediction of Theorem \ref{law:2}.

Finally, in the neighborhood of a BP, we have the following result:
\begin{theorem}[{\rm\cite{HLL2019}}]\label{law:6}
Suppose that $\hat\xi\in (\underline{\xi}, \overline{\xi})$ is a $\mathrm{BP}$ {and \eqref{eq:csss} holds}. If $\overline{\xi}-\underline{\xi}$ is sufficiently small, then   as $\epsilon \to 0$, any positive steady state $u_\epsilon$ of \eqref{eq:10} satisfies
$$
u_{\epsilon}(x,\xi) \to \delta(\xi - \underline\xi)\cdot {u}_1(x) + \delta(\xi - \overline\xi)\cdot {u}_2(x)\quad\text{in distribution},
$$
where $({u}_1, {u}_2)$ is a positive solution of 
\begin{equation}\label{eq:9-a}
\left\{
\aligned
&\underline\xi\partial_{xx} {u}_{1} - \alpha\partial_x{u}_{1}+ {u}_1(r -{u}_1 -{u}_2) = 0 && x\in (0,L),\\
&\overline\xi\partial_{xx} {u}_{2} - \alpha\partial_x{u}_{2} + {u}_2(r -{u}_1 -{u}_2) = 0 && x\in(0,L),\\
&\underline\xi\partial_{x} {u}_{1} - \alpha{u}_1=\overline\xi\partial_{x} {u}_{2} - \alpha{u}_2=0 && x=0,L.
\endaligned
\right.
\end{equation}
\end{theorem}
Theorem \ref{law:6} reveals a new  phenomenon: In the neighborhood of a BP, no single
trait can dominate the competition; instead, the two extreme traits{, $\underline\xi$ and $\overline\xi$,} form a coalition
that together dominates the competition such that any species with other traits cannot invade.
To illustrate Theorem \ref{law:6}, we consider the case 
$r(x)=b_1 e^{c_1 x}+b_2 e^{c_2 x}$ in \eqref{eq:6} with positive constants $b_i$ and  $c_i$, $i=1,2$.
We observe that 
if $\underline\xi=\alpha/c_1$ and $\overline\xi=\alpha/c_2$,
then \eqref{eq:9-a} has the unique positive solution  $u_i(x)=b_i e^{c_i x}$,  $i=1,2$.
Again, in the framework of ideal free distribution, it is  not difficult to understand the biological reasoning: In view of
$u_1(x)+u_2(x)= b_1e^{c_1x}+b_2e^{c_2x}=r(x),$
the total density of the two species can jointly reach an ideal free distribution, so that  $\underline\xi=\alpha/c_1$ and $\overline\xi=\alpha/c_2$ potentially become a pair of ESS. We refer to \cite{CCL2010, CCL2012a, CCL2012b, CCLS2017, Gejji} for related work 
on the evolution of dispersal and ideal free distribution.  

An open problem is whether \eqref{eq:10} or \eqref{eq:9-a} has at most one positive solution, and 
if it exists, whether it is
globally asymptotically stable.
Another challenging question 
is to determine the limit of the time-dependent
solutions $u_\ep(x,\xi,t)$ of \eqref{eq:8-1}
and \eqref{eq:10} 
as $\epsilon \to 0$. See \cite{LLP2020, Nordmann2020}
and references therein for 
recent progress on the rigorous 
derivations of the 
 canonical equations
for the trait evolution
in spatially structured mutation-selection models.


\section{\bf Dynamics of phytoplankton}\label{S6}

Besides the models discussed above, there are  many other types of reaction-diffusion models in  
population dynamics. 
In this section, 
 we briefly discuss some problems arising from the  phytoplankton growth.


Phytoplankton are microscopic plant-like photosynthetic organisms drifting in lakes and oceans and are the foundation of the marine food chain. Since they transport significant amounts of atmospheric carbon dioxide into the deep oceans, they play a crucial role in climate dynamics and have been one of the central topics in marine ecology.
Nutrients and light are the essential resources for the growth of phytoplankton.  
Since most phytoplankton are heavier than water,  
they will sink into the bottom of the lakes or oceans where the light intensity is too weak for their growth. So how do these phytoplankton persist in the water columns? 
Some biologists 
{propose that biased movement, combining with water turbulence, can help phytoplankton diffuse to a position closer to the top of lakes or  oceans, so that they can get access to sunlight.} 
In this connection, 
a series of reaction-diffusion models including  single and multiple phytoplankton species 
are introduced in \cite{HAES, HOW1, HPKS} and the references therein to model the spatio-temporal  dynamics of phytoplankton growth.

The following system of  reaction-diffusion-advection equations was used by Huisman et al. 
\cite{HPKS} to describe the population dynamics of two phytoplankton species:
\begin{equation}\label{equ31}
\left\{
\aligned
&\partial_{t} u=D_1 \partial_{xx} u-\alpha_{1} u_{x}+[g_{1}(I(x, t))-d_{1}]u &&0<x< L,\,\, t >0,\\
&\partial_{t} v=D_2 \partial_{xx} v-\alpha_{2} v_{x}+[g_{2}(I(x, t))-d_{2}]v  &&0<x< L,\,\, t>0,\\
& D_1\partial_{x} u (x, t)-\alpha_{1} u(x, t)=0  &&x=0,L, \,\,t >0, \\
& D_2 \partial_{x} v(x, t)-\alpha_{2} v(x, t)=0 &&x=0,L, \,\,t >0,\\
& u(x,0)=u_{0}(x)\geq, \not \equiv  0,\quad  v(x,0)=v_{0}(x)\geq, \not \equiv 0 && \quad 0\leq x\leq L. 
\endaligned
\right.
\end{equation}
Here $u(x,t), v(x,t)$ denote the population density of the phytoplankton species at depth $x$ and time 
$t$, and  $D_1, D_2> 0 $ are their diffusion rates.  For $i=1,2$,  
$\alpha_{i}\in \mathbb {R} $ is the sinking ($\alpha_{i} >0$) or buoyant ($\alpha_{i} <0$) velocity,  $d_{i}>0 $ is the death rate, and $g_{i}(I) $ represents the specific growth rate of  phytoplankton species as a function of  light intensity $I$,
given by 
the Lambert-Beer law 
\begin{equation}\label{equ2-2}
I(x, t)=I_{0}\exp \left[-k_{0}x-\int_{0}^{x}(k_{1}u(s, t)+k_{2}v(s, t))\mathrm{d}s\right],
\end{equation}
where $I_{0}>0$ is the incident light intensity, $k_{0}> 0$ is the background turbidity that summarizes light absorption by all non-phytoplankton components, and $k_{i}$ is  the absorption coefficient of the $i$-th phytoplankton species. 
Function $g_{i}(I)$ is smooth and satisfies
\begin{equation*}\label{equ3}
g_{i}(0)=0 \quad \text{and} \quad g^{\prime}_{i}(I)>0\,\,\text{ for }\,\, I \geq 0.
\end{equation*}

\subsection{Single phytoplankton species}\label{S5.1}
The single species model, i.e. $v=0$ in \eqref{equ31}, was first considered  in \cite{Shigesada1981} for the self-shading case  and  infinite long water  column; see also
\cite{Ishii1982, Kolokolnikov2009}. 
The authors \cite{DM2010} considered the global dynamics of the single species model when diffusion and drift rates 
{are}  
spatially dependent. 
In \cite{DHL2015}  the authors studied the effect of photoinhibition on the single phytoplankton species, and they found that,  in contrast to the case of no photoinhibition, 
the model with photoinhibition possesses at least two positive steady states in certain parameter ranges. Hsu et al. \cite{Hsu2017} studied the single species growth consuming inorganic carbon with internal storage in a poorly mixed habitat, building upon some interesting recent findings of 
{the principal eigenvalue of a 1-homogeneous positive compact operator.}
In \cite{Peng2015, Peng2016}, the authors considered  the effect of time-periodic
light intensity  at the surface.
Ma and Ou \cite{Ma2017} recently made an important finding that the biomass of the single  species satisfies a comparison principle, even though the density itself 
does not obey such orders.

In \cite{HsuLou2010},  we investigated a single phytoplankton model and obtained some necessary and sufficient conditions for the growth of phytoplankton, and 
the critical death rate, critical water column depth, critical sinking or buoyant coefficient and critical turbulent diffusion rate were studied respectively. One of the results is that 
{the phytoplankton population persists if and only if} the sinking  velocity of  phytoplankton is less than a critical value. There are  many simplified assumptions in the model, e.g.  the death rate  is assumed to be a constant. An interesting issue is to consider the case when the death rate is non-constant in space and time.  Some preliminary results show that this situation is quite different from that of constant death rate.  For instance, under some conditions 
phytoplankton can persist if and only if  the sinking  velocity stays  within an intermediate range, rather than below a single critical value.

\subsection{Two phytoplankton species}\label{S5.1-a}
The coexistence of two or multiple phytoplankton is an important issue.
In contrast to single phytoplakton species, very few results exist for two or more  phytoplankton species; see \cite{Du20081, Du20082, Mei2012}. 
On the one hand, the  coexistence of many species can often be observed  in reality. 
On the other hand, the classical competition theory shows that only the most dominant phytoplankton persists. They seem to contradict each other. The reason is that the classical competition theory is generally established for ordinary differential equation models, i.e. it is assumed that the diffusion rates of phytoplankton are sufficiently  large, so that  only their average densities are considered; i.e. the water column is well mixed.
A natural question is whether  small diffusion can increase the possibility of coexistence.
In this connection, we recently investigated the outcome of 
model \eqref{equ31}-\eqref{equ2-2} in \cite{JLLW2019} and established the following results: 

\smallskip

\begin{itemize}

\item[\rm (i)] If two phytoplankton differ only in their sinking  velocities ($D_1=D_2, \alpha_1\neq \alpha_2, g_1=g_2$, $d_1=d_2$), in \cite{JLLW2019} we showed that the phytoplankton  with smaller sinking  velocity has the competitive advantage. In terms of evolution, the mutation of  phytoplankton can continually  reduce  their sinking  velocities, so that the 
{density} of phytoplankton and water will get closer 
as the population evolves. In other words, as the environment permits, 
the natural selection may favor the phytoplankton whose 
{density} is lighter that of water in some circumstances.

\smallskip

\item[\rm(ii)]
If  two phytoplankton differ only in their diffusion rates ($D_1\neq D_2, \alpha_1=\alpha_2, g_1=g_2$, $d_1=d_2$), it is shown in \cite{JLLW2019} that 
slower diffusion rate will  be selected  when 
buoyant, and in contrast,  faster diffuser wins  when they are sinking with large velocity. An underlying biological reasoning is that when the phytoplankton are buoyant, the slower diffuser are more likely to reach the top of the water column (i.e. without water turbulence),  where the light intensity is the strongest, while when
 sinking with  large velocity, faster diffusion can counterbalance
the tendency to sink and provide individuals with better access to light. 
\end{itemize}

An important tool is the extension of the comparison 
{principle} in \cite{Ma2017} for single-species models to two-species phytoplakton models as follows: 
 \begin{theorem}\label{THm}
Suppose $\{(u_i,v_i)\}_{i=1,2}$ are  non-negative solutions of \eqref{equ31}-\eqref{equ2-2} such that 
\begin{equation*}
\int_0 ^x u_1(s,0)\,\mathrm{d}s
\le, \not\equiv \int_0 ^x u_2(s,0)\,\mathrm{d} s
\quad \text{ and }\quad 
\int_0 ^x v_1(s,0)\,\mathrm{d}s
\ge, \not\equiv \int_0 ^x v_2(s,0)\,\mathrm{d}s
\end{equation*}
{hold}  for all $x\in (0, L]$.
Then for all $x\in (0, L]$ and $t>0$,
\begin{equation*}
\int_0 ^x u_1(s,t)\,\mathrm{d} s
< \int_0 ^x u_2(s,t)\,\mathrm{d} s 
\quad \text{ and }\quad  
\int_0 ^x v_1(s,t)\,\mathrm{d} s
> \int_0 ^x v_2(s,t)\,\mathrm{d} s.
\end{equation*}
\end{theorem}



By Theorem \ref{THm}, system \eqref{equ31}-\eqref{equ2-2} is a strongly monotone dynamical system with respect to  the order generated by  the 
cone {$\mathcal{K}:=\mathcal{K}_1\times (-\mathcal{K}_1)$},
where 
$$
{\mathcal{K}_1}:=\Big\{  \phi \in C([0,L]; \mathbb{R}): \int_{0}^{x} \phi(s)\,\mathrm{d}s\geq 0, \, \forall x\in(0,L] 
\Big\}.
$$
This in turn enables the application of the theory of strongly monotone dynamical system, which  
provides a powerful tool to investigate the global dynamics of two-species system \eqref{equ31}-\eqref{equ2-2}. 

The above findings in 
part (ii) naturally  lead to the following conjecture: 

\begin{conjecture}\label{conjecture4-a} Suppose that $\alpha_1=\alpha_2:=\alpha$, $g_1=g_2$ and $d_1=d_2$.
There exist two critical sinking velocities $\alpha_{\min}$ and  $\alpha_{\max}$ such
that 
the slower diffuser  wins if $ \alpha<\alpha_{\min}$, the faster diffuser  wins if $ \alpha>\alpha_{\max}$, and two species  coexist if $\alpha\in(\alpha_{\min},\alpha_{\max})$.  
\end{conjecture}

\section{\bf Spatial dynamics of epidemic diseases}
\label{S7}

The COVID-19 pandemic  
has impacted or changed 
almost everyone's life. There are many mysteries  about the novel coronavirus,
among which the 
{effect and possible control strategies} of
the movement of individuals is an emerging question. 
The COVID-19 pandemic has been spreading so fast, at least partially  due to the movement of those infected individuals who {are showing little or no symptoms.}
What is the general impact of movement on the persistence and extinction of a disease in spatially heterogeneous environment? 
We studied various susceptible-infected-susceptible (SIS)  models in the heterogeneous environment \cite{Allen2007, Allen2008, Allen2009}, including the following 
reaction-diffusion  model \cite {Allen2008}, where  the susceptible and infected individuals move randomly
and the disease transmission and recovery rates  could be {uneven} 
across the space:
\begin{equation}\label{liu13}
\left\{
\aligned
&\partial_{t}S-d_S\Delta S=-\frac{\beta SI}{S+I}+\gamma I  &&(x,t)\in\Omega\times(0,\infty),\\
&\partial_{t}I-d_I\Delta I=\frac{\beta SI}{S+I}-\gamma I  &&(x,t)\in\Omega\times(0,\infty),\\
&\frac{\partial S}{\partial\nu}=\frac{\partial I}{\partial\nu}=0  &&(x,t)\in\partial\Omega\times(0,\infty),\\
&S(x,0)=S_0(x), \,\,I(x,0)=I_0(x) && x\in \Omega,
\endaligned
\right.
\end{equation}
where  $\Omega$ is a bounded domain in $\mathbb{R}^N$ with smooth boundary $\partial\Omega$ and $\nu$ denotes the
unit outward normal vector on $\partial\Omega$. Here $S(x, t)$ and $I(x, t)$  represent the density of susceptible and infected populations at location $x$ and time $t$, respectively.  Parameters $d_S, d_I>0$ denote their  diffusion rates.
The
functions $\beta(x),\gamma(x)$ account for the rates of disease transmission and recovery,  respectively. They are assumed to be positive and 
at least one of them is non-constant to reflect the spatial heterogeneity of the
residing habitat. 

In \cite {Allen2008}, we defined the basic reproduction number $\mathcal{R}_0$ for the spatial SIS  model \eqref{liu13}, which 
is consistent 
with the next generation approach for heterogeneous populations \cite{Diekmann, van, WZ2012}. 
It is of practical significance  to consider the effect of  spatial heterogeneity of disease
transmission and recovery rates on  the basic reproduction number. For instance, in the study of dengue fever, it is shown in  \cite{WangZhao2011} that  the infection risk may be underestimated if the spatially averaged parameters are used to compute the basic reproduction number for spatially heterogeneous infections.
We proved that  $\mathcal{R}_0$ 
is  monotone decreasing  with respect to the diffusion rate of the infected individuals. 
The basic reproduction number $\mathcal{R}_0$ characterizes the infection risk of disease and serves as
  the threshold value for the extinction and persistence of disease: Namely, if $\mathcal{R}_0<1$, the unique disease-free equilibrium is globally  stable, and if  $\mathcal{R}_0>1$, the disease-free equilibrium is unstable and there is a unique endemic equilibrium. A standing  open question is whether this unique endemic equilibrium is globally stable.

  Interestingly, we also found in \cite{Allen2008} that a disease may persist but it can be controlled  by  limiting  the  movement  of  the  susceptible populations. To be more specific,
  if the spatial environment can be transformed to include low-risk sites
  where the recovery rates are greater than transmission rates (such as through vaccination and treatment), and the movement of the susceptible individuals  can be restricted  (such as through isolation or lock down), then it may be possible to control the disease by flattening the disease outbreak curve, i.e. decreasing 
  the daily number of infected individuals.
 
  In recent years there have been many studies  on SIS and other 
  {types of disease transmission} models
  in spatially and/or temporally varying environments; see \cite{CLL2017, CL2016,  D2019, DW2016,  GD2019,
  GD2020, GR2011, GKLZ2015, 
   LPW2017, LPX2020, MWW2019, Peng2009, PengLiu2009, PengYi2013,  PengZhao2012, WZ2016}. 
   For instance, 
   to study the effect of the movement of exposed individuals on disease outbreaks, 
    the following 
    SEIRS (susceptible-exposed-infected-recovered-susceptible) epidemic reaction-diffusion model 
    was considered in \cite{SLX2019}:
\begin{equation}
\left\{
\aligned
&\partial_t{S}=d_{S}\Delta S-\frac{\beta(x)SI}{S+I+E+R}+\alpha R && (x,t)\in\Omega\times(0,\infty),\\
&{\partial_t{E}}=d_{E}\Delta E+\frac{\beta(x)SI}{S+I+E+R}-\sigma E&& (x,t)\in\Omega\times(0,\infty),\\
&{\partial_t{I}}=d_{I}\Delta I+\sigma E-\gamma(x) I&&(x,t)\in\Omega\times(0,\infty),\\
&{\partial_t{R}}=d_{R}\Delta R+\gamma(x)I-\alpha R && (x,t)\in\Omega\times(0,\infty),\\
&\frac{\partial S}{\partial\nu}=\frac{\partial E}{\partial\nu}=\frac{\partial I}{\partial\nu}=\frac{\partial R}{\partial\nu}=0 
&&(x,t)\in\partial\Omega\times(0,\infty).\\
\endaligned
\right.
\label{e1}
\end{equation}
Here  we divide the individuals into four different
compartments: susceptible ($S$), exposed ($E$), infectious ($I$),  recovered (immune by vaccination, $R$).
The susceptible individuals are infected by infectious individuals with a rate of $\beta$, and become exposed; exposed individuals become infectious with a rate $\sigma$; infected individuals are recovered with a rate $\gamma$; {recovered} individuals lose immunity and go back into the susceptible class with a rate of $\alpha$. Here $S(x,t)$, $E(x,t)$, $I(x,t)$ and $R(x,t)$ denote the density of susceptible, exposed,  infected and recovered individuals at location $x$ and time $t$, and $d_{S}$, $d_{E}$, $d_{I}$, $d_{R}$ represent the  diffusion rates for susceptible, exposed, infected and recovered populations, respectively. 
We assume that the disease transmission rate $\beta(x)$ and recovery rate $\gamma(x)$ are environmentally dependent and could be  spatially heterogeneous.
Our results in \cite{SLX2019} reveal that the travel of exposed individuals could have an important impact on the persistence of disease and the movement of recovered individuals  may enhance the endemic. 
Hence, further 
understanding of the behaviors of the exposed and recovered individuals could be important in designing effective disease control measures.

In \cite{LL2020} we investigated the spatial SIS model \eqref{liu13} with spatially heterogeneous and time-periodic coefficients 
and proved that  the basic reproduction number $\mathcal{R}_0$ is non-decreasing with respect to the time-period $T$.
In some scenario, this would imply that there is a critical period $T^*>0$ such that $\mathcal{R}_0<1$ for $T<T^*$ and  $\mathcal{R}_0>1$ for $T>T^*$; i.e. increasing the period of disease transmissions and recovery  will increase the chance of disease outbreak. This might suggest why
the outbreak of some epidemic diseases are less frequent than others, e.g. annually vs bi-annually, as different epidemic diseases could have different transmission and recovery rates so that the threshold values $T^*$ are different.

How do movement and spatial heterogeneity affect the competition among multiple strains?
Many diseases are {present in} multiple phenotypes or strains, e.g. the 
novel coronavirus has found three major 
strains 
which are 
prevalent in different continents. 
It is an important evolutionary problem to determine the environmental conditions under which the adaptability of a certain strain  becomes dominant.
Bremermann and Thieme \cite{Bremermann} showed that for an 
SIR model with one host and multiple pathogens,  different pathogens will die out eventually except those that optimize the basic reproduction number. 
Recent studies  showed that for some non-autonomous multi-strain models, even if the transmission and recovery rates change  time-periodically, 
these strains cannot coexist; i.e. the temporal heterogeneity may not increase the chance of the coexistence for multiple pathogens. Tuncer and Martcheva \cite{Tuncer2012} {and Wu et al. \cite{WTM2017} considered} the following  two-strain SIS spatial model  to address the question  whether the presence of spatial structure would allow two strains to coexist, as the corresponding spatially homogeneous model generally leads to competitive exclusion:  \begin{equation}\label{liu13-11}
\left\{
\aligned
&\partial_{t}S-d_S\Delta S=-\frac{S(\beta_1 I_1+\beta_2 I_2)}{S+I_1+I_2}+\gamma_1 I_1+\gamma_2I_2 &&(x,t)\in\Omega\times(0,\infty),\\
&\partial_{t}I_i-d_{I,i}\Delta I_i=\frac{\beta_i SI_i}{S+I_1+I_2}-\gamma_i I_i &&(x,t)\in\Omega\times(0,\infty),
\\
&\frac{\partial S}{\partial\nu}=\frac{\partial I_i}{\partial\nu}=0  &&(x,t)\in\partial\Omega\times(0,\infty),\\
&S(x,0)=S_0(x), \,\,I_i(x,0)=I_{0,i}(x) && x\in \Omega,
\endaligned
\right.
\end{equation}
where $i=1,2$. These two strains could be different in dispersal rates, transmission rates and recovery rates.
For the two  pathogens differing only in their diffusion rates, {a conjecture is that} the pathogen with the slower diffusion will drive the other one  to extinction for proper transmission and recovery rates; {see \cite{WTM2017} for some partial results}. For the two  pathogens differing only in their  recovery rates, we conjecture that the strain with the recovery rate 
of greater spatial variation will be selected eventually.
The coexistence of multiple pathogens in spatially heterogeneous 
and temporally varying environment remains a promising open research direction. 
\bigskip
\bigskip

\noindent{\bf Acknowledgment.} 
We thank the anonymous referee and Gr\'{e}goire Nadin and Lei Zhang for their helpful comments. KYL and YL  are partially supported by the NSF grant DMS-1853561;
SL was partially supported by the Outstanding Innovative Talents Cultivation Funded Programs 2018 of Renmin Univertity of China and the NSFC grant No. 11571364. 
\bigskip

\baselineskip 18pt
\renewcommand{\baselinestretch}{1.2}

\end{document}